\documentclass[review]{elsarticle}

\usepackage{amsmath}
\usepackage{amsfonts}
\usepackage{amssymb}
\usepackage{amsthm}
\usepackage{bm}
\usepackage{indentfirst}
\usepackage{graphicx}
\usepackage{subfigure}
\usepackage{array}
\usepackage{fullpage}

\DeclareMathAlphabet{\mathsfsl}{OT1}{cmss}{m}{sl}

\newcommand{\PreserveBackslash}[1]{\let\temp=\\#1\let\\=\temp}
\newcolumntype{C}[1]{>{\PreserveBackslash\centering}p{#1}}
\newcolumntype{R}[1]{>{\PreserveBackslash\raggedleft}p{#1}}
\newcolumntype{L}[1]{>{\PreserveBackslash\raggedright}p{#1}}

\numberwithin{equation}{section}
\newtheorem{thm}{Theorem}[section]
\newtheorem{lem}[thm]{Lemma}
\newtheorem{cor}[thm]{Corollary}
\theoremstyle{definition}

\newtheorem{rem}[thm]{Remark}
\newtheorem{assu}[thm]{Assumption}

\usepackage{anyfontsize}
\usepackage{enumerate}
\usepackage{esint}
\usepackage{etoolbox}
\usepackage{isomath}
\usepackage{mathrsfs}
\usepackage{mathtools}
\usepackage{multicol}
\usepackage{pgfplots}
\usepackage{scalerel}
\usepackage{stackengine}
\usepackage{tabulary}
\usepackage{tikz}

\makeatletter
\newcommand*\bdot{\mathpalette\bdot@{.65}}
\newcommand*\bdot@[2]{\mathbin{\vcenter{\hbox{\scalebox{#2}{$\m@th#1\bullet$}}}}}
\makeatother

\makeatletter
\newcommand*\bddot{\mathpalette\bddot@{.65}}
\newcommand*\bddot@[2]{\mathbin{\vcenter{\hbox{\scalebox{#2}
    {$\m@th#1\smash{{}_{\bullet}^{\bullet}}$}}}}}
\makeatother

\newcommand{\circled}[2][]{%
  \tikz[baseline=(char.base)]{%
    \node[shape = circle, draw, inner sep = .5pt]
    (char) {\phantom{\ifblank{#1}{#2}{#1}}};%
    \node at (char.center) {\makebox[0pt][c]{#2}};}}
\robustify{\circled}

\makeatletter
\newcommand{\opnorm}{\@ifstar\@opnorms\@opnorm}
\newcommand{\@opnorms}[1]{%
  \left|\mkern-1.5mu\left|\mkern-1.5mu\left|
   #1
  \right|\mkern-1.5mu\right|\mkern-1.5mu\right|
}
\newcommand{\@opnorm}[2][]{%
  \mathopen{#1|\mkern-1.5mu#1|\mkern-1.5mu#1|}
  #2
  \mathclose{#1|\mkern-1.5mu#1|\mkern-1.5mu#1|}
}
\makeatother

\stackMath
\newcommand\reallywidecheck[1]{%
\savestack{\tmpbox}{\stretchto{%
  \scaleto{%
    \scalerel*[\widthof{\ensuremath{#1}}]{\kern-.6pt\bigwedge\kern-.6pt}%
    {\rule[-\textheight/2]{1ex}{\textheight}}
  }{\textheight}%
}{0.5ex}}%
\stackon[1pt]{#1}{\scalebox{-1}{\tmpbox}}%
}
\usepackage[utf8]{inputenc}
\usepackage{bm}
\usepackage{mathrsfs}
\usepackage{amsmath}
\usepackage{amsfonts}
\usepackage{amsthm}
\usepackage{algorithm}
\usepackage[noend]{algpseudocode}
\usepackage{color}
\usepackage{setspace}
\usepackage{float}
\usepackage{import}
\usepackage{url}
\usepackage{multicol}
\usepackage{subfigure}
\biboptions{sort&compress}
\usepackage[top=1in,bottom=1in,left=1in,right=1in]{geometry}
\newcommand{\real}{\mathbb{R}}

\newcommand{\R}{\mathbb{R}}
\newcommand{\C}{\mathbb{C}}

\newcommand{\mcD}{\mathcal{D}}

\newcommand{\mcG}{\mathcal{G}}

\newcommand{\mcL}{\mathcal{L}}

\newcommand{\mcA}{\mathcal{A}}

\usepackage{xcolor}

\def\omg{{\Omega}}

\def \bb{\bm{b}}

\def \thetab{{\boldsymbol \theta}}
\def \betab{{\boldsymbol\beta}}

\def \fb{\bm{f}}

\def \ub{\bm{u}}

\def \vb{\bm{v}}
\def \xb{\bm{x}}

\def \zb{\bm{z}}

\def \yb{\bm{y}}

\def \eb{\bm{e}}

\def \xib{{\boldsymbol\xi}}

\newcommand{\vertii}[1]{{\left\vert\left\vert #1
    \right\vert\right\vert}}    
    
\newcommand{\verti}[1]{{\left\vert #1
    \right\vert}}  
    
\begin{document}

\begin{frontmatter}

\title{An asymptotically compatible probabilistic collocation method for randomly heterogeneous nonlocal problems}

\address[yy]{Department of Mathematics, Lehigh University, Bethlehem, PA, 18015}
\address[xt]{Department of Mathematics, University of California, San Diego, CA, 92093}
\address[xy]{Department of Industrial and Systems Engineering, Lehigh University, Bethlehem, PA, 18015}
\address[xl]{Department of Mathematics and Statistics, University of North Carolina at Charlotte, Charlotte, NC, 28223}
\address[cw]{Department of Mathematics, University of Texas at Austin, Austin, TX, 78712}

\author[yy]{Yiming Fan}\ead{yif319@lehigh.edu}
\author[xt]{Xiaochuan Tian}\ead{xctian@ucsd.edu}
\author[xy]{Xiu Yang}\ead{xiy518@lehigh.edu}
\author[xl]{Xingjie Li}\ead{xli47@uncc.edu}
\author[cw]{Clayton Webster}\ead{claytongwebster@utexas.edu}
\author[yy]{Yue Yu}\ead{yuy214@lehigh.edu}


\begin{abstract}
In this paper we present an asymptotically compatible meshfree method for solving nonlocal equations with random coefficients, describing diffusion in heterogeneous media. In particular, the random diffusivity coefficient is described by a finite-dimensional random variable or a truncated combination of random variables with the Karhunen-Lo\`{e}ve decomposition, then a probabilistic collocation method (PCM) with sparse grids is employed to sample the stochastic process. On each sample, the deterministic nonlocal diffusion problem is discretized with an optimization-based meshfree quadrature rule. We present rigorous analysis for the proposed scheme and demonstrate convergence for a number of benchmark problems, showing that it sustains the asymptotic compatibility spatially and achieves an algebraic or sub-exponential convergence rate in the random coefficients space as the number of collocation points grows. Finally, to validate the applicability of this approach we consider a randomly heterogeneous nonlocal problem with a given spatial correlation structure, demonstrating that the proposed PCM approach achieves substantial speed-up compared to conventional Monte Carlo simulations.
\end{abstract}


\end{frontmatter}


\tableofcontents

\section{Introduction}

Since the last decade, there has been a great interest in using nonlocal integro-differential equations to describe physical systems, due to their natural ability to describe physical phenomena at small scales and their reduced regularity requirements which lead to greater flexibility \cite{silling_2000,bazant2002nonlocal,zimmermann2005continuum,emmrich2007analysis,
emmrich2007well,zhou2010mathematical, du2011mathematical,du2016multiscale,podlubny1998fractional,mainardi2010fractional,magin2006fractional,
burch2011classical,du2014nonlocal,defterli2015fractional,lischke2018fractional,du2014peridynamics,antoine2005approximation,dayal2007real,sachs2013priori,chiarello2018global,erbay2018convergence,bucur2016nonlocal,you2020data}. These nonlocal models are defined in terms of a lengthscale $\delta$, referred to as a horizon, which denotes the extent of nonlocal interaction. The nonlocal viewpoint allows a natural description of processes requiring reduced regularity in the relevant solution, such as the peridynamics model for fracture mechanics \cite{bazant2002nonlocal,du2013nonlocal,yu2018partitioned}. An important feature of such models is that they revert back to corresponding classical partial differential equation (PDE) models as the horizon $\delta \rightarrow 0$. When refining the spatial discretization characterized by grid size $h$ such that $h \rightarrow 0$, discretization methods which preserve correct local limits are termed asymptotically compatible (AC) schemes \cite{tian2014asymptotically}, and there has been significant work in recent years toward establishing such discretizations, see, e.g.,  \cite{tian2014asymptotically,d2020numerical,leng2021asymptotically,pasetto2018reproducing,hillman2020generalized,seleson2016convergence,du2016local,trask2019asymptotically,You_2019,you2020asymptotically,tao2017nonlocal}. Broadly, strategies either involve adopting traditional weak form via finite element shape functions and carefully performing geometric calculations to integrate over relevant horizon/element subdomains, or adopting a strong form meshfree discretization where particles are associated with abstract measure \cite{yu2021asymptotically}. The former is more amenable to mathematical analysis due to a better variational setting, while the latter is simple to implement and generally with smaller computational cost \cite{silling2005meshfree,bessa2014meshfree}. In this work we pursue the asymptotically compatible meshfree approach.

One of the limitations of the current {state-of-art} works is that most of them consider a homogenized nonlocal model, which may not work well when the material is heterogeneous and its microsctructure plays a critical role. In a recent study \cite{zhao2020stochastic}, Zhao et al. found that a fully homogenized peridynamic model fails to capture certain correct fracture modes/patterns in reinforced concrete. Therefore, they have proposed a stochastic bond-based peridynamic model where the material property is described as random fields. {The type of each bond connecting material points $\xb$ and $\yb$ was modeled by a random variable, and the discrete probability distribution of this random variable depends on the volume fraction of aggregate and cement on $\xb$ and $\yb$.} With this model, fracture patterns and the order in which various cracks develop match experimental observations. Their findings indicate the importance of considering the spatial variability of material properties in nonlocal models, especially when the physical parameters describing spatially varying properties of heterogeneous materials cannot be accurately characterized in all details.

In the present paper, we consider a stochastic nonlocal diffusion equation, where the heterogeneous material property is modeled by a random field. The solution of this stochastic equation describes the probability density function (PDF) of the state variable, e.g., the concentration. Differing from \cite{zhao2020stochastic} which studied the solution pattern on each individual realization rather than the solution statistics, we focus on the numerical estimates of the first two statistical moments, i.e., the mean and (co)variance. The mean provides an unbiased estimate of the variables and the variance quantifies the uncertainty associated with this estimate. The Monte Carlo (MC) method and its variations \cite{stein1987large,loh1996latin,fox1999strategies,niederreiter1992random,cliffe2011multilevel} are usually used to solve stochastic moment equations and often considered a reliable numerical tool \cite{dagan1998comment}. In the MC method, a large number (denoted by $K$) of random realizations (samples) are generated for the prescribed random inputs and repetitive deterministic solvers are employed for each sample. The results are then statistically analyzed based on all $K$ samples to calculate leading moments of variables of interest. However, its slow $O(K^{-1/2})$ convergence rate hinders the application of MC method on relatively large scale problems, since one has to solve the differential equation for every sample. Moreover, comparing with local (classical) PDE models, numerically solving nonlocal equations is often more expensive due to its relative lack of sparsity. Therefore, the need for efficient and accurate stochastic numerical methods is even more pressing in the nonlocal setting.

To achieve a faster convergence rate, several stochastic numerical methods were developed for stochastic local (classical) PDE models, including probabilistic Galerkin methods (PGMs) \cite{babuska2004galerkin,babuvska2005solving,ghanem2003stochastic,le2004uncertainty,matthies2005galerkin,xiu2002wiener,wan2005adaptive}, probablistic collocation methods (PCMs) \cite{xiu2005high,nobile2008anisotropic,ma2009adaptive,zhang2012error,lin2009efficient}, reduced basis methods \cite{rozza2007reduced,rozza2007stability,chen2014comparison,chen2013weighted,elman2013reduced,guan2017reduced}, etc. Among these methods, the probabilistic collocation method with sparse grids inherits the ease of implementation in the MC methods since only solutions at sample points are needed. At the same time, it also reduces the required number of sample points to achieve a given numerical accuracy, especially on problems with small random dimensions and  {sufficient solution smoothness in the parameter space}. Therefore, in this work we will employ the probabilistic collocation methods with sparse grids. Comparing with the attentions received by stochastic PDE problems, numerical studies of nonlocal problems in the uncertainty quantification setting remain limited. In \cite{guan2017reduced}, reduced-basis methods are developed for constructing surrogates of the solution of a parameterized nonlocal diffusion problem with random input data in a finite element method framework. However, to the authors' best knowledge, there exists no work on studying the solution smoothness and the theoretical limiting behavior of stochastic nonlocal problems when $\delta\rightarrow 0$, while these studies are crucial to the design of accurate and asymptotically compatible stochastic numerical schemes. Moreover, the application and rigorous error estimates of meshfree discretization method also remain limited for stochastic nonlocal problems.

The major contribution of the present work is to propose a complete workflow of asymptotically compatible stochastic numerical methods and rigorous mathematical analysis for randomly heterogeneous nonlocal diffusion problems. In particular, we propose to employ a meshfree method with optimization-based quadrature rule \cite{trask2019asymptotically,yu2021asymptotically} for the discretization in the physical space, and a PCM with sparse grids for the discretization in the parameter space. By proving that {the solution of the nonlocal equation with diffusivity coefficient described as a finite-dimensional random field is analytic in the input random variables}, we show that the sparse grid PCM achieves at least algebraic convergence with the increase of sample points. Moreover, given sufficiently large level of sparse grid formulation, the sparse grid PCM converges sub-exponentially. To characterize the convergence in physical space and the asymptotic compatibility, we for the first time provide analysis for the stochastic nonlocal diffusion problem with random coefficients, showing that its solution converges to the local solution when $\delta\rightarrow 0$. Based on these mathematical analysis, we provide error estimates for the optimization-based quadrature rule with both fixed horizon $\delta$ and also $\delta$ going to an asymptotic limit. Lastly, we develop a complete workflow to solve for randomly heterogeneous nonlocal  problems, by representing the heterogeneous material coefficient as a random field with given spatial correlation structure and approximating the coefficient by a truncated combination of random variables using the Karhunen-Lo\`{e}ve expansion. This work provides a road map to add uncertainty quantification functionality onto pre-existing asymptotically compatible code for deterministic nonlocal problems in a non-intrusive manner, which achieves algebraic or sub-exponential convergence in the solution mean and variance while sustaining the spatial asymptotic compatibility to the correct local limit.

We remark that the paper is organized to establish the rigorous mathematical underpinnings of the approach in the first half, while the second half focuses on a numerical verification and more engineering-oriented exploration of its application. The paper is organized as follows, with all major notations listed in Table \ref{tab:notation}. We recall first the relevant results in nonlocal calculus and provide mathematical analysis for the deterministic and stochastic nonlocal diffusion problems in Section \ref{sec:pde}. After establishing the continuous limits of the stochastic nonlocal problem, we next pursue a consistent discretization. In Section \ref{sec:method}, we propose our numerical approach for stochastic nonlocal problems by employing the sparse grid PCM and an optimization-based meshfree quadrature rule in the physical space and establish rigorous error estimates. In particular, we show that the proposed approach achieves algebraic or sub-exponential convergence in the parametric space, and address the convergence rates to the nonlocal and local limits respectively. The theoretical error estimates are verified on a number of one-dimensional and two-dimensional problems with analytic solutions for the local and nonlocal limits in Section \ref{sec:veri}. In Section \ref{sec:exp}, we further extend the proposed formulation to handle a more engineering-oriented problem, where the random diffusivity coefficient is modeled by a random field with a given spatial correlation structure. Section \ref{sec:conclusion} summarizes our findings and discusses future research.

\section{The Deterministic and Parametric Nonlocal Diffusion Problem}\label{sec:pde}

\begin{table}[]
    \centering
    \begin{tabular}{|c|l|}
    \hline
    Symbol&Description\\
    \hline
    $\omg$&Physical domain.\\
    $\delta$&Horizon size.\\
    $B_\delta(\xb)=\{\yb:|\yb-\xb|\leq\delta\}$&The physical interaction region surrounding $\xb$.\\
    $\omg_\delta$&Nonlocal boundary which is a collar of thickness surrounding $\omg$.\\
    $(\Omega_p,\mathcal{F},\mathcal{P})$&Probability space.\\
    $\Gamma$&Space of random variables.\\
    $d$&Dimension of the physical space $\omg$.\\
    $N$&Dimension of the random space, i.e., the number of random variables.\\
    $M$&Total number of grid points for spatial discretization.\\
    $K$&Total number of collocation points (samples) in PCM.\\
    $\mathcal{D}$&Nonlocal divergence operator.\\
    $\mathcal{G}$&Nonlocal gradient operator.\\
    $\mcL^\delta$&Nonlocal diffusion operator.\\
    $\mcL^0$&Local (classical) diffusion operator.\\
    $\mcL_h^\delta$&Discretized nonlocal diffusion operator.\\
    $v_{(i)}$& The $i$-th component of vector $\vb$.\\
    $M_{(i,j)}$& The $i$-th row $j$-th column element of matrix $\bm{M}$.\\
    $u^D:\omg_\delta\rightarrow \real$&Dirichlet type boundary condition.\\
    $\gamma(\xb,\yb)=\gamma_\delta(|\yb-\xb|)$&Symmetric kernel function.\\
    $s$&Order of singularity in kernel $\gamma$.\\
    $A:(\omg\cup\omg_\delta)\times(\omg\cup\omg_\delta)\times\Gamma\rightarrow\real$&Nonlocal random diffusion strength function.\\
    $a:\omg\times\Gamma\rightarrow\real$&Local random diffusion strength function.\\
    $r$ and $R$&Lower and upper bounds of $A$ and $a$.\\
    $S_\delta(\omg)$ and $T_{\delta}$&Nonlocal energy space and the corresponding bilinear form.\\
    $\rho:\Gamma\rightarrow\real^+$&Probability density of the random variable.\\
    $\chi_h=\{\xb_i\}$ and $h$&(Quasi-uniform) grid set and the grid size for spatial discretization.\\
    $\bm{V}_h$&Space of functions to be integrated exactly in the spatial discretization.\\
    $\omega_{j,i}$&Quadrature weight for $\xb_j$ to generate integral in $B_\delta(\xb_i)$.\\
    $\Theta_N=\{\xib_k\}$&Prescribed nodes for the Lagrange interpolation in the random space $\Gamma$.\\
    $\mu_k$&Corresponding quadrature weight for $\xib_k$ in the random space.\\
    $\zeta$&Sparseness parameter for the Smolyak sparse grid formulation.\\
    $\eta=\zeta-N$&Level of the Smolyak sparse grid formulation.\\
    $\phi^\epsilon$& Mollification function.\\
    $u^\delta:(\omg\cup\omg_\delta)\times\Gamma\rightarrow \real$&Nonlocal solution of \eqref{eqn:nonlocal_random}.\\
    $u^0:\omg\times\Gamma\rightarrow \real$&Local solution of \eqref{eqn:local_random}.\\
    $u^\delta_h$&Numerical solution for the deterministic nonlocal diffusion problem.\\
    $u_{h,K}^\delta$&Numerical solution with spatial grid set $\chi_h$ and sparse grid level $\eta$ in PCM.\\
    $\Xi(\xb,\yb)$&Covariance kernel for the variances between material points $\xb$ and $\yb$.\\
    \hline
    \end{tabular}
    \caption{Table of Notations.}
    \label{tab:notation}
\end{table}

In this section we introduce the major notations and definitions will be used throughout this paper. We begin with Table \ref{tab:notation} and in Section \ref{sec:diff} we introduce the deterministic nonlocal diffusion problem while Section \ref{sec:randomdiff} is dedicated to the stochastic nonlocal diffusion problem. Moreover, we provide novel theoretical analysis for the stochastic nonlocal diffusion problem, namely, its compatibility with the classical companion and analytic regularity.

\subsection{Nonlocal Calculus and Deterministic Nonlocal Diffusion Problem}\label{sec:diff}

In this section, we review the the governing equations of deterministic nonlocal diffusion models which provide the foundation for the stochastic nonlocal problems of interest. Given that $\Omega\subset \mathbb{R}^d$, $d\in \mathbb{Z}^+$, is a bounded Lipschitz domain, we consider the nonlocal elliptic equation in $\Omega$. To do so, we first introduce the relevant nonlocal calculus. Let $\bm{\alpha}(\xb,\yb):\real^d\times\real^d\rightarrow \real^d$ be an antisymmetric function, for a vector function  $\bm{v}(\xb,\yb):\real^d\times\real^d\rightarrow \real^d$, we define the nonlocal divergence $\mathcal{D}[\bm{v}]:\real^d\rightarrow\real$:
$$\mathcal{D}[\bm{v}](\xb):=\int_{\real^d}\left(\bm{v}(\xb,\yb)+\bm{v}(\yb,\xb)\right)\cdot\bm{\alpha}(\xb,\yb)d\yb, \quad \xb\in\real^d,$$
and for a scalar function $u(\xb):\real^d\rightarrow\real$ we define the nonlocal gradient $\mathcal{G}[u]:\real^d\times\real^d\rightarrow\real^d$:
$$\mathcal{G}[u](\xb,\yb):=\left(u(\yb)-u(\xb)\right)\bm{\alpha}(\xb,\yb), \quad \xb,\yb\in\real^d.$$
As shown in \cite{du2013nonlocal}, the adjoint operator of $\mathcal{D}$ {with respect to the $L^2$ inner product} is $\mathcal{D}^*=-\mathcal{G}$. We then consider a nonlocal diffusion problem where every point $\xb\in\omg$ is interacting with a neighborhood of points, and their interaction is described by a symmetric kernel function $\gamma(\xb,\yb):=\bm{\alpha}(\xb,\yb)\cdot \bm{\alpha}(\xb,\yb)$ and a two-point scalar function $A(\xb,\yb)$ representing the nonlocal diffusion strength. {Without loss of generality, we assume in this paper that $A$ is symmetric in its two arguments, i.e., $A(\xb,\yb)= A(\yb,\xb)$. 
} A nonlocal diffusion operator on a scalar function $u:\real^d\rightarrow \real$ is then given by 
$$\mcL[u](\xb):=\mathcal{D}[A(\xb,\yb)\mathcal{G}[u]](\xb)=2\int_{\real^d}A(\xb,\yb)\gamma(\xb,\yb)(u(\yb)-u(\xb))d\yb, \quad \xb\in\real^d.$$
In this paper we further assume that the interacting kernel function $\gamma$ is radial and compactly supported on a Euclidean ball surrounding $\xb$, i.e., $B_\delta(\xb):=\{\yb\in\real^d:|\yb-\xb|<\delta\}$:
\begin{equation}\label{eqn:require_ga}
\left\{\begin{array}{l}
\gamma(\xb,\yb)=\gamma_\delta(|\xb-\yb|)=\frac{1}{\delta^{d+2}}\gamma_1\left(\frac{|\xb-\yb|}{\delta}\right)=\frac{D_0}{\delta^{d+2-s}|\xb-\yb|^s}\\
\text{ where $\gamma_1$ is a nonnegative and nonincreasing function with $s$-th order singularity, satisfying}\\
\text{supp}(\gamma_1)\subset B_1(\bm{0}) \text{ and }\int_{B_1(\bm{0})}\gamma_1(|\zb|)|\zb|^2d\zb=d.
\end{array}\right.
\end{equation}
The above kernel assumptions have implications on the boundary conditions that are prescribed on a collar of thickness $\delta$ outside the domain $\omg$, that we denote as 
$$\omg_\delta:=\left\{ \bm x\in \real^d\backslash\omg:  \text{dist}(\bm x, \partial\omg)  <\delta\right\}$$
and refer to as nonlocal boundary. Therefore, the static nonlocal diffusion problem in a deterministic parameters setting is given as:
\begin{equation}\label{eqn:nonlocal}
\left\{\begin{array}{ll}
-\mcL^\delta[u](\xb)=-\mathcal{D}[A(\xb,\yb)\mathcal{G}[u]](\xb)=f(\xb),\quad &\text{ for }\xb\in\omg\\
u(\xb)=u^D(\xb),\quad &\text{ for }\xb\in\omg_\delta,
\end{array}\right.
\end{equation}
where $u^D$ is the given Dirichlet-type boundary datum in the nonlocal trace space \cite{Du2021trace}. Without loss of generality, for the analysis, we consider homogeneous Dirichlet boundary conditions $u^D(\xb)=0$, and the proposed method is applied to inhomogeneous Dirichlet-type problems in numerical tests of Section \ref{sec:veri}. Note that although the proposed model can be applied to other boundary conditions, e.g., the Neumann-type boundary conditions in \cite{You_2019,you2020asymptotically}, here we focus on the Dirichlet-type nonlocal constraint problem for simplicity.

We assume that the diffusion coefficient $A\in L^\infty((\omg\cup\omg_\delta)^2)$. To make the nonlocal diffusion problem \eqref{eqn:nonlocal} uniformly elliptic, we make the uniform boundedness assumption on the diffusion coefficient $A(\xb,\yb)$: 
\begin{equation}\label{eqn:require_A}
0<r\leq A(\xb,\yb)\leq R<\infty,\quad \text{ for }\xb,\yb\in\omg\cup\omg_\delta,
\end{equation}
where $r$ and $R$ are two positive constants. 
\eqref{eqn:nonlocal} is then associated with the nonlocal energy semi-norm
$$| u |^2_{S_\delta(\omg)}= \int_{\omg\cup\omg_\delta}  \int_{\omg\cup\omg_\delta} \gamma_\delta(| \bm y-\bm x|) (u(\bm y)-u(\bm x))^2 d\bm y d\bm x$$
where the nonlocal energy space is defined as
$$S_\delta(\omg)=\left\{u\in L^2(\omg\cup\omg_\delta):\int_{\omg\cup\omg_\delta}  \int_{\omg\cup\omg_\delta} \gamma_\delta(| \bm y-\bm x|) (u(\bm y)-u(\bm x))^2 d\bm y d\bm x<\infty,\,u|_{\omg_\delta}=0\right\}.$$
Moreover, we define the  bilinear form $T_\delta:S_\delta\times S_\delta\rightarrow\real$ as
$$T_\delta[v,w]=\int_{\omg\cup\omg_\delta}  \int_{\omg\cup\omg_\delta} A(\bm x,\bm y) \gamma_\delta(|\xb-\yb|)(v(\bm y)-v(\bm x))(w(\bm y)-w(\bm x)) d\bm y d\bm x.$$
With the boundedness assumption proposed in \eqref{eqn:require_A} 
we have $0<r\verti{v}^2_{S_\delta(\omg)}\leq T_\delta[v,v]\leq R\verti{v}^2_{S_\delta(\omg)}$. Notice that to solve for the weak solution of \eqref{eqn:nonlocal}, we find $u(\xb)\in S_\delta(\omg)$ such that
\begin{equation}\label{eqn:weak_nonlocal}
T_\delta[u,v]=(f,v)_{L^2(\omg)},\quad \forall v\in S_\delta(\Omega).
\end{equation}
As proved in \cite{mengesha2014bond}, the bilinear form $T_\delta$ holds uniform ellipticity as well as the nonlocal Poincar\'{e} inequality, and therefore the weak formulation for the deterministic nonlocal equation \eqref{eqn:weak_nonlocal} is well-posed.
Here we summarize the relevant results in the following lemma:
\begin{lem}\label{Lemma:wellposedness_nonlocal}{{\cite[Proposition~5.3]{mengesha2014bond}}}
Assume that $\gamma(\xb,\yb)$ satisfies the conditions in \eqref{eqn:require_ga}, then there exist generic constants $C$ and $\delta_0>0$ such that for all $0<\delta<\delta_0$, the nonlocal energy semi-norm $|\cdot|_{S_\delta(\omg)}$ satisfies the nonlocal Poincar\'{e} inequality 
\begin{equation}
\vertii{v}^2_{L^2(\omg)}\leq C |v|^2_{S_\delta(\omg)}, \quad \forall v\in S_\delta(\omg). 
\end{equation}
\end{lem}

As a result of the Poincar\'e inequality, the semi-norm $|\cdot|_{S_\delta(\omg)}$ defined is also a norm on $S_\delta(\omg)$. For the rest of the paper, we will use $\| \cdot\|_{S_\delta(\omg)}:=| \cdot|_{S_\delta(\omg)}$ to denote the norm on $S_\delta(\omg)$. 
Combining Lemma \ref{Lemma:wellposedness_nonlocal} with the properties of $A$ in \eqref{eqn:require_A}, we can see that the bilinear form $T_\delta$ is bounded and coercive:
\begin{equation}
T_\delta[v,w]\leq C\vertii{v}_{S_\delta(\omg)}\vertii{w}_{S_\delta(\omg)}, \quad \forall v,w\in S_\delta(\omg),    
\end{equation}
\begin{equation}
T_\delta[v,v]\geq C\vertii{v}^2_{S_\delta(\omg)}, \quad \forall v\in S_\delta(\omg).   
\end{equation}
Therefore, by the Lax-Milgram theorem, there exists a unique solution $u^\delta\in S_\delta(\omg)$ for the deterministic nonlocal diffusion problem \eqref{eqn:nonlocal} for each $f\in (S_\delta(\omg))^\ast$, where $(S_\delta(\omg))^\ast$ is the dual space of $S_\delta(\omg)$ equipped with the induced norm $\| f\|_{(S_\delta(\omg))^\ast}:=\underset{v\in S_\delta(\omg), v\neq 0}{\sup} \frac{\langle f, v\rangle}{\| v\|_{S_\delta(\omg)}}$.  $\langle \cdot,\cdot\rangle$ denotes the duality pairing between $(S_\delta(\omg))^\ast$ and $S_\delta(\omg)$, and $\langle f, v\rangle = (f,v)_{L^2(\omg)}$ when $f\in L^2(\omg)\subset (S_\delta(\omg))^\ast$.

Next, we consider the compatibility of the nonlocal diffusion and the classical companion. To properly define the local limit of \eqref{eqn:nonlocal} as $\delta\to0$, we need to make the following continuity assumption on the diffusion coefficient $A(\xb,\yb)$: 
\begin{equation}\label{eqn:require_B}
A(\xb,\yb) \in C(\overline{(\omg\cup\omg_\delta)^2}) \text{ and }   a(\xb):=A(\xb,\xb).  
\end{equation}
Therefore, the nonlocal diffusion operator $\mcL^\delta[u]$ has a companion of the classical diffusion operator $-\nabla\cdot(a(\xb)\nabla(u(\xb)))$, and \eqref{eqn:nonlocal} can be seen as a nonlocal analogue to the local diffusion equation with Dirichlet-type boundary condition:
\begin{equation}\label{eqn:local}
\left\{\begin{array}{ll}
\mcL^0 [u](\xb):= -\nabla\cdot(a(\xb)\nabla(u(\xb)))=f(\xb),\quad &\text{ for }\xb\in\omg\\
u(\xb)=u^D(\xb),\quad &\text{ for }\xb\in\partial\omg,
\end{array}\right.
\end{equation}
From \eqref{eqn:require_A}, we have $0<r\leq a(\xb)\leq R$ for any $\xb\in\omg$, so \eqref{eqn:local} has a unique and bounded solution in {$H^1(\omg):=\{u\in L^2(\omg)\big|
\int_\omg |\nabla u(\xb)|^2 d\xb <\infty
\}$ with corresponding boundary condition $u|_{\partial \omg}(\xb)=u^{D}(x)$. When we consider $u^D = 0$, the solution space is denoted by $H_0^1(\omg):=\{ u\in H^1(\omg), u|_{\partial\omg} = 0\}$.}

Notice that we take the minimal regularity assumptions on $A(\xb, \yb)$ and $a(\xb)$.    
Denoting the solution of local problem \eqref{eqn:local} as $u^0(\xb)$ and the solution of the nonlocal one \eqref{eqn:nonlocal} with a given horizon size $\delta$ as $u^\delta(\xb)$, we now show the convergence theorem with the minimal regularity assumptions. 

\begin{thm}
\label{thm:compatibility}
Assume that $\gamma(\xb,\yb)$ satisfies the conditions in \eqref{eqn:require_ga}, and  $A(\xb,\yb)$ satisfies the conditions in \eqref{eqn:require_A}. 
Let $f\in (S_\delta(\omg))^\ast$, the dual space of $S_\delta(\omg)$, then we have 
\begin{equation}
\label{eqn:energy_estimate}
 \| u^\delta\|_{S_\delta(\omg)}\leq  \frac{\| f\|_{(S_\delta(\omg))^\ast}}{r}.   
\end{equation}
In addition,  if $\| f\|_{(S_\delta(\omg))^\ast}$ is uniformly bounded for all $\delta\in (0,\delta_0)$ and $A(\xb,\yb)$ satisfies \eqref{eqn:require_B},  then the nonlocal and local diffusion problems are compatitble as $\delta \rightarrow 0$:
$$\underset{\delta\rightarrow 0}{\lim}\vertii{u^\delta-u^0}_{L^2(\omg)}=0.$$
\end{thm}
\begin{proof}
We first show the proof of \eqref{eqn:energy_estimate}. Since $u^\delta$ is a solution to the nonlocal problem, we have
\[
T_\delta[u^\delta, v] = \langle f,v\rangle \leq \| f\|_{(S_\delta(\omg))^\ast} \| v\|_{S_\delta(\omg)}
\]
for any test function $v\in S_\delta(\omg)$. Now let $v= u^\delta$, we get $r\| u^\delta\|^2_{S_\delta(\omg)}\leq T_\delta[u^\delta, u^\delta] \leq \| f\|_{(S_\delta(\omg))^\ast} \| u^\delta\|_{S_\delta(\omg)}$. 
Therefore, we have \eqref{eqn:energy_estimate}. 

The proof of the second part involves two steps. In the first step, we assume that $A(\xb, \yb)\in C^\infty( (\omg\cup\omg_\delta)^2)$. 
Then it is easy to see by Taylor expansion that for any $v\in C_0^\infty(\omg)$ (with zero extended values on $\omg_\delta$), we have the pointwise convergence of   $\mathcal{D}[A(\xb,\yb)\mathcal{G}[u]](\xb) $ to $\nabla \cdot(a \nabla v(\xb) )$ as $\delta\to0$. Indeed,  by doing Taylor expansion of $A$ and $v$ around $\xb$, we find that 
\[
\begin{split}
&\mcD[A \mcG [v]](\xb)= 2\int A(\xb, \yb)  \gamma_\delta(|\yb-\xb|) (v(\yb)-v(\xb)) d\yb \\
=& 2 \int \Big(a(\xb) + (\yb-\xb) \cdot \frac{\nabla a(\xb)}{2} \Big)  \gamma_\delta(|\yb-\xb|) \Big( (\yb-\xb)^T \nabla v(\xb)  + 
\frac{1}{2}(\yb-\xb)^T D^2 v(\xb) (\yb-\xb)\Big) d\yb +O(\delta^2)\\
=&  \nabla a(\xb)\cdot \nabla v(\xb)  + a(\xb)\Delta v(\xb) +O(\delta^2)\,.
\end{split}
\]
Note that the above equalities are obtained since  $\int_{B_\delta(\xb)}\gamma_\delta(|\xb-\yb|)\underset{i_1+\cdots+i_d=3}{\Pi}(\yb-\xb)^{i_k}_{(k)}d\yb=0$ thanks to the symmetry of the kernel $\gamma_\delta$, where $(\yb-\xb)_{(k)}$ denotes the $k$-th component of $(\yb-\xb)$ and $i_k\geq 0$ is the power on that component.
Then we argue that the convergence is also in $L^2(\omg)$ since $|\mcD(A \mcG v )(\xb) |$ is uniformly bounded for $v\in C_0^\infty(\omg)$ and $\delta\in (0,\delta_0)$.  Notice that from the assumption on $\| f\|_{(S_\delta(\omg))^\ast}$, we have  $\| u^\delta\|_{S_\delta(\omg)}$ being uniformly bounded for all $\delta\in (0,\delta_0)$. Then using similar arguments in \cite{tian2014asymptotically}, we can show $\|u^\delta -u^0 \|_{L^2(\omg)}\to 0$ as $\delta\to0$.

 For the general case that $A\in C(\overline{(\omg\cup\omg_\delta)^2})$, we will use the mollification technique.
 Take standard mollifiers $\phi^\epsilon \in C^\infty(\R^{2d})$, we define $A^\epsilon = \phi^\epsilon \ast A$.  We denote the solution to \eqref{eqn:nonlocal} associated with coefficient $A^\epsilon$ to be $u^{\delta, \epsilon}$. Then we can use the first step to conclude that 
$
\|u^{\delta, \epsilon} - u^{0, \epsilon} \|_{L^2(\omg)} \xrightarrow{\delta\to0} 0
$, where $u^{0,\epsilon}$ is the solution to \eqref{eqn:local} associated with coefficient $a^\epsilon(\xb) : = A^\epsilon(\xb,\xb)$. Now in order to show $\| u^\delta - u^0\|_{L^2(\omg)}\to 0$, we notice that
\[
\lim_{\delta\to0}\| u^\delta - u^0\|_{L^2(\omg)} \leq \sup_{\delta\in (0,\delta_0)} \| u^{\delta, \epsilon}-u^{\delta} \|_{L^2(\omg)} +  \lim_{\delta\to0}\|u^{\delta, \epsilon} - u^{0, \epsilon} \|_{L^2(\omg)} +  \| u^{0, \epsilon}-u^{0} \|_{L^2(\omg)} ,
\]
for any $\epsilon>0$.
Therefore, we only need to show 
\begin{equation}\label{eqn:step2estimates}
\left\{ 
\begin{aligned}
& \lim_{\epsilon\to0}\sup_{\delta\in (0,\delta_0)} \| u^{\delta, \epsilon}-u^{\delta} \|_{L^2(\omg)}=0, \quad \text{ and } \\
&\lim_{\epsilon\to0} \| u^{0, \epsilon}-u^{0} \|_{L^2(\omg)}=0. 
\end{aligned} 
\right.
\end{equation}
For the first equation in \eqref{eqn:step2estimates}, we first notice that $\| A^\epsilon -A\|_{C((\omg\cup\omg_\delta)^2)} \to 0 $ as $\epsilon\to0$ since $A(\xb, \yb)$ is uniformly continuous on $\overline{(\omg\cup\omg_\delta)^2}$. Now since $u^{\delta, \epsilon}$ and $u^\delta$ are solutions to \eqref{eqn:nonlocal} with different coefficients and the same right-hand side, we have 
\[
(A^\epsilon \mcG [u^{\delta, \epsilon} - u^\delta], \mcG [v])_{L^2((\omg\cup\omg_\delta)^2)} = ((A-A^\epsilon) \mcG [u^{\delta, \epsilon}], \mcG [v])_{L^2((\omg\cup\omg_\delta)^2)} = : \langle g^{\delta,\epsilon}, v\rangle , 
\]
for any $v\in S_\delta(\omg)$. We can show $\langle g^{\delta,\epsilon}, v \rangle \to 0$ as $\epsilon\to0$ uniformly independent of $\delta$ since 
\[
\langle g^\epsilon, v\rangle  \leq \| A-A^\epsilon\|_{C((\omg\cup\omg_\delta)^2)} \| u^{\delta, \epsilon} \|_{S_\delta(\omg)} \| v \|_{S_\delta(\omg)}  \leq C  \| A-A^\epsilon\|_{C((\omg\cup\omg_\delta)^2)}  \| v \|_{S_\delta(\omg)} ,
\]
where we have used $ \| u^{\delta, \epsilon} \|_{S_\delta(\omg)} \leq C$ from \eqref{eqn:energy_estimate}. 
Therefore, we have 
\[
\sup_{\delta\in (0,\delta_0)}\| u^{\delta,\epsilon} - u^\delta \|_{S_\delta(\omg)} \leq C\sup_{\delta\in (0,\delta_0)} \| g^{\delta, \epsilon}\|_{(S_\delta(\omg))^\ast} \leq C \| A-A^\epsilon\|_{C((\omg\cup\omg_\delta)^2)} \to 0
\] 
as $\epsilon\to0$ and the convergence in $L^2$ is then implied from the Poincar\'e inequality. 
The proof for the second equation in \eqref{eqn:step2estimates} is similar by noticing that $\| a^{\epsilon}- a\|_{C(\omg)}\to0$ as $\epsilon\to0$. 
\end{proof}

\subsection{Parametric Nonlocal Diffusion Problem}\label{sec:randomdiff}

We now consider the case in which the coefficient $A$ is provided by a random field  $A(\xb,\yb,\omega)$, where $\omega\in\Omega_p$ and $\Omega_p$ is the sample space of a probability space $(\Omega_p, \mathcal{F}, \mathcal{P})$. Here,
$\mathcal{F}$ is the $\sigma$-algebra of subsets of $\Omega_p$ and $\mathcal{P}$ is the probability measure. In practice, this random field is usually represented in a ``truncated'' form using a limited number of random variables (see an example in Section~\ref{sec:KL}). Thus, it can be rewritten as $A(\xb, \yb, \xib)$, where $\xib=(\xi_{(1)},\xi_{(2)},\dotsc, \xi_{(N)})$, $N$ is a positive integer {which denotes the dimension of the parametric space}, and $\xi_{(i)}$ are random variables. In practice, we often assume they are independent and identically distributed (i.i.d.) random variables.
Under this setting, we consider $A(\xb,\yb,\xib):(\omg\cup\omg_\delta)\times (\omg\cup\omg_\delta)\times\Gamma\rightarrow \real$, where $\Gamma$ is the space of $\xib$ and it is typically called random space or parametric space. 
Without loss of generality, here we assume that $\Gamma= \prod_{i=1}^N \Gamma_{i} \subset \R^N$ where $\Gamma_{i} = [-1,1]$, and the random variable $\xib\in \Gamma$ has a probability density $\rho:\Gamma \to \R^+$. We are interested in solving the family of nonlocal elliptic equations given by 
\begin{equation}\label{eqn:nonlocal_random}
\left\{\begin{array}{ll}
-\mathcal{D}[A(\xb,\yb,\xib)\mathcal{G}[u]](\xb)=f(\xb),\quad &\text{ for }\xb\in\omg\\
u(\xb)=u^D(\xb),\quad &\text{ for }\xb\in\omg_\delta.
\end{array}\right.
\end{equation}
For each $\xib\in \Gamma$, we assume that $A(\xb, \yb, \xib)$ is symmetric in its first two variables and $A(\xb, \yb, \xib) \in L^\infty((\omg\cup\omg_\delta)^2)$.
We also assume the uniform ellipticity of the nonlocal problems, i.e., 
\begin{equation}\label{eqn:UniformElliptic}
    0<r\leq A(\xb,\yb,\xib)\leq R<\infty. 
\end{equation}
Therefore the Lax-Milgram theorem ensures the well-posedness of nonlocal problem for each $\xib\in \Gamma$. 
In addition, in order to consider the limit $\delta\to0$, we need to assume that for each $\xib\in\Gamma$, 
\begin{equation}\label{eqn:DiagCont}
A(\cdot,\cdot, \xib)\in C(\overline{(\omg\cup\omg_\delta)^2}) \text{ and }   a(\xb,\xib):=A(\xb,\xb,\xib) .   
\end{equation}
Then we have the corresponding family of local elliptic equation 
\begin{equation}\label{eqn:local_random}
\left\{\begin{array}{ll}
 -\nabla\cdot\big(a(\xb,\xib)\nabla(u(\xb))\big)=f(\xb),\quad &\text{ for }\xb\in\omg\\
u(\xb)=u^D(\xb),\quad &\text{ for }\xb\in\partial\omg.
\end{array}\right.
\end{equation}

For each given parameter $\xib \in \Gamma$, we denote the solution to the nonlocal equation \eqref{eqn:nonlocal_random} by $u^\delta(\xb, \xib)$ and the solution to the corresponding local equation \eqref{eqn:local_random} by  $u^0(\xb, \xib)$. A corollary of Theorem \ref{thm:compatibility} is that $u^\delta(\xb, \xib)$ converges to $u^0(\xb, \xib)$ in the space $L^2(\omg)\otimes L^2_{\rho}(
\Gamma)$ as $\delta\to0$. 
\begin{cor}
Assume that $A(\xb,\yb,\xib)$ satisfies \eqref{eqn:UniformElliptic} and \eqref{eqn:DiagCont}, then we have 
{\[
\lim_{\delta\rightarrow 0}\| u^\delta - u^0\|_{L^2( \omg)\otimes  L^2_\rho(\Gamma)}= 0\,.
\]}
\end{cor}
\begin{proof}
For any $\xib\in \Gamma$, and, we know from Theorem \ref{thm:compatibility} that  $\| u^\delta(\cdot, \xib)\|_{S_\delta(\omg)}\leq C $ for all $\delta\in (0,\delta_0)$ and  $\| u^\delta(\cdot, \xib) - u^0(\cdot, \xib)\|_{L^2(\omg)}\to 0$ as $\delta\to0$. Therefore, it is easy to see that $\| u^\delta(\cdot, \xib) - u^0(\cdot, \xib)\|_{L^2(\omg)} \leq C$ for all $\xib\in\Gamma$ and $\delta\in (0,\delta_0)$.  By invoking the dominated convergence theorem, we have
\[
\| u^\delta - u^0 \|_{ L^2(\omg) \otimes L^2_\rho(\Gamma) } 
= \int_{\Gamma} \| u^\delta(\cdot , \xib)-u^0( \cdot, \xib) \|^2_{L^2(\omg)}   \rho(\xib)d\xib \overset{\delta\rightarrow 0}{\longrightarrow} 0 \,.
\]
\end{proof}

To discuss the regularity of solutions with respect to the parameter space, we need to assume the existence of a holomorphic extension of $A(\xb, \yb,\xib)$. 
\begin{assu}[Holomorphic parameter dependence]
\label{assu:continuation}
The complex continuation of $A(\xb,\yb, \xib)$, represented as the map $A: \C^N \to L^\infty((\omg\cup\omg_\delta)^2)$, is a $L^\infty((\omg\cup\omg_\delta)^2)$-valued holomorphic
function on $\C^N$.
\end{assu}

This condition is easily fulfilled with $A(\xb, \yb, \xib)$ consisting of polynomials, exponential,
sine and cosine functions of the variables $\xi_{(1)}, \xi_{(2)}, \cdots, \xi_{(N)}$. For example, the holomorphic extension exists if  $A(\xb,\yb, \xib) = \bar A(\xb,\yb) +\sum_{i=1}^N f_i(\xi_{(i)}) \psi_i(\xb,\yb)$, where $\bar A \in L^\infty((\omg\cup\omg_\delta)^2)$, $\psi_i \in L^\infty((\omg\cup\omg_\delta)^2)$, and $f_i$ is a polynomial, exponential, sine or cosine function ($1\leq i\leq N$).   

\subsubsection{Analytic regularity}
In order for the function $u^\delta(\xb, \xib)$ to be uniformly recovered by polynomial expansions in the parameter space, we will show the analyticity of the solution
$u^\delta$ with respect to the parameterization. By Assumption \ref{assu:continuation}, we can extend the definition of $A(\xb,\yb,\xib)$ to $A(\xb,\yb,\hat\xib)$ where $\hat\xib$ belongs to the complex domain
\[
\hat\Gamma := \otimes_{1\leq i\leq N} \{ \hat\xi_{(i)} \in \C : |\hat\xi_{(i)}|\leq 1 \}. 
\]
Next, we need the complex uniform ellipticity assumption, namely there exists  $r_c, R_c\in \R$ such that 
\begin{equation}\label{eqn:ComplexUniformElliptic}
    0<r_c\leq\text{Re}(A(\xb,\yb,\hat\xib)) \leq |A(\xb,\yb,\hat\xib)|\leq R_c<\infty. 
\end{equation}
for $\hat\xib\in \hat\Gamma$.
Therefore the nonlocal problem with the complex coefficient $A(\xb,\yb,\hat\xib)$ is well-posed and the corresponding solution is denoted by $u^\delta(\xb, \hat\xib)$. Here we also remark that for the analyticity of solutions to hold, non-affine coefficients may also be used as long as $\xib\mapsto A(\xb,\yb,\xib)$ possesses an analytic extension to the complex domains and the complex uniform ellipticity condition is satisfied, see related discussion for the local elliptic equations in \cite{tran2017analysis}. We also define the set 
\[
\mcA : = \{\hat\xib \in \C^N : \frac{r_c}{2} < \text{Re}(A(\xb,\yb,\hat\xib)) \leq |A(\xb,\yb,\hat\xib)|< 2R_c \}. 
\]
So it is clear that the set $\mcA$ contains $\hat\Gamma$. 

For the discussions in this Section, we need the complex function spaces. We let $S_\delta(\omg, \C)$ be the space of complex valued functions with norm
\[
\| u \|^2_{S_\delta(\omg,\C)} =  \int_{\omg\cup\omg_\delta}  \int_{\omg\cup\omg_\delta} \gamma_\delta(| \bm y-\bm x|) |u(\bm y)-u(\bm x)|^2 d\bm y d\bm x, 
\]
where $|u(\bm y)-u(\bm x)|^2$ is understood as $\overline{(u(\bm y)-u(\bm x))} (u(\bm y)-u(\bm x))$ and $\gamma_\delta$ is still the real valued kernel function.
$(S_\delta(\omg, \C))^\ast$ is the dual space of  $S_\delta(\omg,\C)$ equipped with the induced norm. The main result of the this Section is the following theorem on the analyticity of solutions. 

\begin{thm}\label{thm:analytic}
Assume that $A(\xb,\yb,\hat\xib)$ satisfies \eqref{eqn:ComplexUniformElliptic}, then the function $\hat\xib \mapsto u(\xb,\hat\xib)$ is holomorphic in an open neighborhood of $\hat\Gamma$. 
 \end{thm}
 
 To prove the theorem, we first need a stability lemma. 
 \begin{lem} Suppose $u^\delta$ and $\tilde u^\delta$ are two solutions of \eqref{eqn:nonlocal} with the same right-hand side $f$ and different coefficients $A(\xb,\yb)$ and $\tilde A(\xb,\yb)$ satisfying \eqref{eqn:require_A}, then  
\[
\| u^\delta-  \tilde u^\delta \|_{S_\delta(\omg)} \leq   \| A-\tilde A\|_{L^\infty((\omg\cup\omg_\delta)^2)} \frac{\| f \|_{(S_\delta(\omg))^\ast}}{r^2}  
\]
\end{lem}
\begin{proof}
Subtracting the variational formulations for $u^\delta$ and $\tilde u^\delta$, we find that
\[
\begin{split}
0 = \iint_{(\omg\cup\omg_\delta)^2} \Big(A \mcG [u^\delta](\xb,\yb)- \tilde A \mcG [\tilde u^\delta]&(\xb,\yb)\Big)\mcG [v](\xb,\yb) d\yb d\xb \\
 =  \iint_{(\omg\cup\omg_\delta)^2}  A \big( \mcG [u^\delta](\xb,\yb)- \mcG [\tilde u^\delta](\xb,\yb)\big) &\mcG [v](\xb,\yb) d\yb d\xb  +
 \iint_{(\omg\cup\omg_\delta)^2} (A- \tilde A) \mcG [\tilde u^\delta](\xb,\yb) \mcG [v](\xb,\yb) d\yb d\xb \,.
\end{split}
\]
Therefore $w= u^\delta- \tilde u^\delta$ is a solution of 
\[
 \iint_{(\omg\cup\omg_\delta)^2} A(\xb,\yb) \mcG[w](\xb,\yb) \mcG [v](\xb,\yb) d\yb d\xb = \langle l, v \rangle 
\]
where $\langle l, v \rangle := \iint ( \tilde A- A)   \mcG[\tilde u^\delta] \mcG[v] d\yb d\xb$. So from \eqref{eqn:energy_estimate} we have 
\[
\| w\|_{S_\delta(\omg)}\leq \frac{\| l\|_{(S_\delta(\omg))^\ast}}{r} \leq \frac{\| A-\tilde A\|_{L^\infty((\omg\cup\omg_\delta)^2)} \| \tilde u^\delta\|_{S_\delta(\omg)}}{r}\leq \| A-\tilde A\|_{L^\infty((\omg\cup\omg_\delta)^2)} \frac{\| f \|_{(S_\delta(\omg))^\ast}}{r^2}.
\]
\end{proof}
From the proof of the stability lemma, it is easy to see that similar stability estimate holds when $A$ and $\tilde A$ are dependent on the parameter $\xib\in \R^N$ (or $\hat\xib\in \C^N$) as long as $A$ (or Re$A$) and $\tilde A$ (or Re$\tilde A$) are bounded below.  
Now Theorem \ref{thm:analytic} can be proved in a similar way as the its local counterpart in 
\cite[Lemma 2.2]{cohen2011analytic}.   

\noindent {\it Proof of Theorem \ref{thm:analytic}.} 
For simplicity of notations, sometimes we only write the explicit dependence of functions on the parameter $\hat\xib$ in the discussion, although they may also depend on the spatial variables. 
First, it is easy to see that $\mcA$ is an open neighborhood of $\hat\Gamma$ by the continuity of the maps $\hat\xib \mapsto \text{Re} A(\hat\xib)$ and $\hat\xib\mapsto |A(\hat\xib)|$. Now we will show that for all $\hat\xib\in \mcA$, 
  the function $\hat\xib \mapsto  u^\delta(\hat\xib)$ admits a complex derivative $\partial_{\hat \xi_{(i)}}  u^\delta(\hat\xib) \in S_\delta(\omg, \C)$ for $i\in \{1,\cdots, N\}$.
  
Let $\eb_i \in \R^N$ be the unit vector in the $i$-th dimension.  For $\hat h\in\C\backslash\{0\}$, define the difference quotient function
\[
w^\delta_{\hat h}(\hat\xib) =\frac{ u^\delta(\hat\xib+\hat h \eb_i)- u^\delta(\hat\xib)}{\hat h} \in S_\delta(\omg;\C)\,.
\]
Since the maps $\hat\xib \mapsto \text{Re} A(\hat\xib)$ and $\hat\xib\mapsto |A(\hat\xib)|$ are continuous, we have the boundedness of $\text{Re} A(\hat\xib+\hat h \eb_i)$ and $|A(\hat\xib+\hat h \eb_i)|$ for  sufficiently small $\hat h$. Therefore $w^\delta_{\hat h}$ is well-defined for sufficiently small $\hat h$. Then for all $ v\in S_\delta(\omg;\C)$, 
\[
\begin{split}
0=&\iint_{(\omg\cup\omg_\delta)^2}  A(\xb,\yb ,\hat\xib+\hat h \eb_j) \mcG[u^\delta(\hat\xib+\hat h \eb_j)](\xb,\yb) \overline{\mcG [ v](\xb,\yb)} d\yb d\xb \\
& -\iint_{(\omg\cup\omg_\delta)^2}  A(\xb,\yb ,\hat\xib) \mcG[ u^\delta(\hat\xib)](\xb,\yb) \overline{\mcG [ v](\xb,\yb)} d\yb d\xb \\
=& \hat h \iint_{(\omg\cup\omg_\delta)^2}  A(\xb,\yb ,\hat\xib)  \mcG[w^\delta_{\hat h}(\hat\xib)](\xb,\yb) \overline{\mcG [ v](\xb,\yb)} d\yb d\xb \\
& +\hat h\iint_{(\omg\cup\omg_\delta)^2} \psi_i(\xb,\yb)  \mcG[ u^\delta(\hat\xib+\hat h \eb_j)](\xb,\yb) \overline{\mcG [ v](\xb,\yb)} d\yb d\xb\,.
\end{split}
\]  
Therefore $w^\delta_{\hat h}$ is the unique solution to the variational problem 
\[
\iint_{(\omg\cup\omg_\delta)^2}  A(\xb,\yb ,\hat\xib)  \mcG[w^\delta_{\hat h}(\hat\xib)](\xb,\yb) \overline{\mcG [ v](\xb,\yb)} d\yb d\xb =\langle l_{\hat h},v\rangle,
\]
where  $\langle l_{\hat h},v\rangle:=\iint_{(\omg\cup\omg_\delta)^2} \psi_i(\xb,\yb)  \mcG[ u^\delta(\hat\xib+\hat h \eb_j)](\xb,\yb) \overline{\mcG [ v](\xb,\yb)} d\yb d\xb $. One can show that $l_{\hat h}$ converges to $l_0$ in $(S_\delta(\omg;\C))^\ast$. Indeed, for all $v\in S_\delta(\omg;\C)$, 
\[
\begin{split}
|\langle l_{\hat h} - l_0, v\rangle| &= \left| \iint_{(\omg\cup\omg_\delta)^2} \psi_i(\xb,\yb)  \mcG[ u^\delta(\hat\xib+\hat h \eb_i)-u^\delta(\hat\xib)](\xb,\yb) \overline{\mcG [ v](\xb,\yb)} d\yb d\xb \right| \\
&\leq \| \psi_i\|_{L^\infty((\omg\cup\omg_\delta)^2)} \| u^\delta(\hat\xib+\hat h \eb_i)-u^\delta(\hat\xib)\|_{S_\delta(\omg;\C)} \| v\|_{S_\delta(\omg;\C)}  \\
&\leq \hat h \| \psi_i\|_{L^\infty((\omg\cup\omg_\delta)^2)} \frac{\| f\|_{(S_\delta(\omg;\C))^\ast}}{(r_c/2)^2} \| v\|_{S_\delta(\omg;\C)} ,
\end{split}
\]
in which the last inequality comes from the stability estimate. Therefore, $w^\delta_{\hat h}$ converges in $S_\delta(\omg;\C)$ to $w^\delta_0$, which is the solution to 
\[
\iint_{(\omg\cup\omg_\delta)^2}  A(\xb,\yb ,\hat\xib)  \mcG[w^\delta_0(\hat\xib)](\xb,\yb) \overline{\mcG [ v](\xb,\yb)} d\yb d\xb =\langle l_0,v\rangle.
\]
Hence $\partial_{\hat \xi_{(i)}}  u^\delta(\hat\xib)=w^\delta_0(\hat\xib) \in S_\delta(\omg;\C)$. 
\qed

\section{Spatial and Stochastic Numerical Methods}\label{sec:method}

\begin{algorithm}
\caption{Overall algorithm for the stochastic nonlocal diffusion problem \eqref{eqn:nonlocal_random}}
{{\begin{algorithmic}[1]
\State{{\bf Offline Stage}: 

$1a)$ For spatial discretization, determine a meshfree grid set $\chi_h=\{\xb_i\}_{i=1}^M\in \omg\cup\omg_\delta$ and calculate the optimization-based quadrature weights $\{\omega_{j,i}\}_j$ for each $\xb_i\in\chi_h$ by solving \eqref{eq:quadQP}.

$1b)$ For probabilistic collocation method, determine the (sparse) collocation points set in the parametric space $\Theta_N=\{\xib_k\}_{k=1}^K\in \Gamma$ and the corresponding quadrature weights $\mu_k$ following \eqref{eqn:sparsesample}.
}
\State{{{\bf Online Stage}: Solving for \eqref{eqn:nonlocal_random}}: For $k=1,\cdots,K$, do

$2a)$ Assemble the stiffness matrix $\bm{Q}=[Q_{(i,j)}]$ such that
\begin{displaymath}
Q_{ij}=\left\{\begin{array}{ll}
2\underset{\xb_j \in \chi_h\cap B_\delta(\xb_i)}{\sum}A(\xb_i,\xb_j,\xib_k)\gamma_\delta(|\xb_i-\xb_j|)\omega_{j,i},& \text{ if }i=j,\\
-2A(\xb_i,\xb_j,\xib_k)\gamma_\delta(|\xb_i-\xb_j|)\omega_{j,i},& \text{ if }i\neq j \text{ and }\xb_j\in B_\delta(\xb_i),\\
0,& \text{ else },
\end{array}\right.
\end{displaymath}
and the right-hand-side vector $\fb=[f(\xb_i)]_{i=1}^M$.

$2b)$ Compute the numerical solution $\ub=[u_h^\delta(\xb_i,\xib_k)]_{i=1}^M$ for the deterministic problem corresponding to the $k$-th sample, by solving $\ub=\bm{Q}^{-1}\fb$.

}
\State{{\bf Postprocessing Stage}: Generate statistical moments of the random solution following \eqref{eqn:E}-\eqref{eqn:std}.

}
\end{algorithmic}}}
\label{algorithm:1}
\end{algorithm}

\subsection{Spatial: Optimization-Based Meshfree Quadrature Rules}\label{sec:meshfree}


In this section we introduce a strong form of particle discretizations of the deterministic nonlocal diffusion problem introduced in Section \ref{sec:diff}. This approach is based upon the optimization-based quadrature rule developed in \cite{trask2019asymptotically,yu2021asymptotically}.  Denoting the numerical solution of \eqref{eqn:nonlocal} as $u^\delta_h$, two types of convergence are desired in the developed numerical scheme:
\begin{equation}\label{eqn:convh}
\underset{h\rightarrow 0}{\lim}\vertii{u_h^\delta-u^\delta}_{L^2(\omg)}=0,\quad \text{ and }\quad\underset{h,\delta\rightarrow 0}{\lim}\vertii{u_h^\delta-u^0}_{L^2(\omg)}=0.
\end{equation}
The first type of convergence indicates that the numerical discretization method is consistent with the nonlocal problem, while the second type shows that the nonlocal numerical solution preserves the correct local limit, or equivalently, the numerical scheme is asymptotically compatible. To maintain an easily scalable implementation, in asymptotic compatibility studies we assume $\delta$ to be chosen such that the ratio $\frac{\delta}{h}$ is bound by a constant as $\delta \rightarrow 0$, restricting ourselves to the ``$\delta$-convergence'' scenario \cite{bobaru2009convergence}\footnote{In some literature it is also denoted as the ``M-convergence'', see, e.g., \cite{yu2021asymptotically}.}. In this setting, one obtains banded stiffness matrices allowing scalable implementations. Typically in the literature a scheme is termed as asymptotically compatible (AC) if it recovers the classical solution for any $\delta,h\rightarrow 0$. Here we abuse the definition slightly and only require the $\delta$-convergence.
We will provide truncation error estimates for the quadrature error convergence rates to the nonlocal analytical solution and the local limit, respectively.

\subsubsection{Mathematical Formulation and Implementation}

Discretizing the whole interaction region $\Omega\cup\omg_\delta$ by a collection of points $\chi_{h} = \{\xb_i\}_{\{i=1,2,\cdots,M\}} \subset \Omega\cup\omg_\delta$, we aim to solve for the solution $u_{(i)}\approx u(\xb_i)$ on all $\xb_i\in \chi_h$. 
Although the method can be applied to more general grids, for analysis we require $\chi_h$ to be a uniform Cartesian grid:
$$\chi_h:=\{(k_{(1)}h,\cdots,k_{(d)}h)|\bm{k}=(k_{(i)},\cdots,k_{(d)})\in\mathbb{Z}^d\}\cap (\omg\cup\omg_\delta).$$
Here $h$ is the spatial grid size. 
For the deterministic nonlocal diffusion model \eqref{eqn:nonlocal} we pursue a discretization through the following one point quadrature rule at $\chi_h$ \cite{silling_2010}:
\begin{align}
    &\nonumber-(\mathcal{L}^\delta_h u)_{(i)}:=-2\sum_{\xb_j \in \chi_h\cap B_\delta(\xb_i)\backslash\{\xb_i\}}A(\xb_i,\xb_j)\gamma_\delta(|\xb_i-\xb_j|)(u_{(j)}-u_{(i)})\omega_{j,i} = f(\xb_i),&\;\text{ for }\xb_i\in\chi_h\cap\omg,\\
    &u_{(i)}=u^D(\xb_i),&\;\text{ for }\xb_i\in\chi_h\cap\omg_\delta,\label{eq:discreteNonloc}
\end{align}
where we specify {$\left\{\omega_{j,i}\right\}$} as a to-be-determined collection of quadrature weights admitting interpretation as a measure associated with each collocation point $\xb_i$. Note that although we only solve for $u_{(i)}$ on grid points in $\chi_h$, we will denote the numerical solution $u_h^\delta$ of the above nonlocal problem as the piecewise constant interpolation of $u_{(i)}$.

We use an optimization-based approach to define these weights extending previous work \cite{trask2019asymptotically,yu2021asymptotically}, constructed to ensure consistency guarantees. Specifically, we seek quadrature weights for integrals supported on balls of the form
\begin{equation}
I[q] := \int_{B_\delta (\xb_i)} q(\xb,\yb) d\yb \approx I_h[q] := \sum_{\xb_j \in \chi_h\cap B_\delta(\xb_i)\backslash\{\xb_i\}} q(\xb_i,\xb_j) \omega_{j,i}
\end{equation}
where we include the subscript $i$ in $\left\{\omega_{j,i}\right\}$ to denote that we seek a different family of quadrature weights for different subdomains $B_\delta(\xb_i)$. Denoting $\bm{P}_m(\real^d)$ as the space of $m$-th order polynomials, we obtain these weights from the following optimization problem
\begin{align}\label{eq:quadQP}
  \underset{\left\{\omega_{j,i}\right\}}{\text{argmin}} \sum_{\xb_j \in \chi_h\cap B_\delta(\xb_i)\backslash\{\xb_i\}} \frac{\omega_{j,i}^2W(|\xb_j-\xb_i|)}{2} \quad
  \text{such that}, \quad
  I_h[q] = I[q] \quad \forall q \in \bm{V}_h,
\end{align}
where $\bm{V}_h= \bm{S}_{\gamma_\delta,\xb}:=\left\{ q = p(\yb-\xb)\gamma_\delta(|\xb-\yb|) \,|\, p(\xb) \in \bm{P}_3(\mathbb{R}^d)\right\}$ denotes the space of functions which should be integrated exactly. Note that when $p(\xb)=const$, we have $q(\xb,\yb)=\gamma_\delta(\verti{\xb-\yb})$. Therefore, reproducing this function requires the kernel function $\gamma_\delta$ to be integrable, or equivalently, $s<d$. {$W(r)$ is a radially symmetric positive weight function supported in $B_\delta(\bm{0})$. Following the discussions in \cite{leng2021asymptotically}, we take $W(r)=\gamma_\delta(r)$.} As shown in \cite{trask2019asymptotically}, this particular choice of reproducing space provides the minimal reproducing set to achieve the optimal $O(\delta^2)$ asymptotic convergence rate in nonlocal problems with homogeneous diffusion coefficient. 
Moreover, we note that provided the quadrature points are unisolvent over the desired reproducing space, \eqref{eq:quadQP} may be proven to have a solution by interpreting it as a generalized moving least squares (GMLS) problem \cite{leng2021asymptotically}. For certain choices of $\bm{V}_h$, such as $m$-th order polynomials, unisolvency holds under the following assumptions: the domain $\Omega$ satisfies a cone condition, the pointset $\chi_h \cap B_\delta(\xb_i)$ is quasi-uniform, and $h/\delta$ is sufficiently small \cite{wendland2004scattered}.

For each $\xb_i\in\chi_h\cap\omg$, we denote the set of its interacting neighbor points as $\{\xb_j\}_{j=1,\cdots,M_i}=\chi_h\cap B_\delta(\xb_i)\backslash\{\xb_i\}$ and the corresponding quadrature weights in a size $M_i$ vector, denoted as $\bm{\omega}=[\omega_{1,i},\cdots,\omega_{M_i,i}]^T$. With the Lagrange multiplier method, quadrature weights may be obtained from \eqref{eq:quadQP} by solving the saddle-point problem for each $\xb_i$
\begin{equation}\label{eqn:QPsaddle}
\begin{bmatrix}
\bm{W} & \bm{H}^\intercal \\
\bm{H} & \bm{0}
\end{bmatrix}
\begin{bmatrix}
\bm{\omega} \\
\bm{\lambda}
\end{bmatrix}
 =
 \begin{bmatrix}
 \bm{0} \\
 \bm{g}
 \end{bmatrix},
\end{equation}
where $\bm{W}$ is an $M_i\times M_i$ diagonal matrix with $\bm{W}_{kk}=\gamma_\delta(|\xb_k-\xb_i|)$, $\bm{\lambda} \in \mathbb{R}^{dim(\bm{V}_h)}$ is a set of Lagrange multipliers used to enforce reproducability, $\bm{H}=[H_{(\alpha,j)}] \in \mathbb{R}^{M_i \times dim(\bm{V}_h)}$ consists of the reproducing set evaluated at each quadrature point (i.e. $H_{(\alpha ,j)} = p_\alpha(\xb_j), \, \text{for all } p_\alpha \in \bm{V}_h$), and $\bm{g}=[g_{(\alpha)}] \in \mathbb{R}^{dim(\bm{V}_h)}$ consists of the integral of each function in the reproducing set over the ball such that $g_{(\alpha)} = I[p_\alpha]$. In particular, when considering $\bm{V}_h= \bm{S}_{\gamma_\delta,\xb}$, and denoting $\gamma_{ji}:=\gamma_\delta(|\xb_i-\xb_j|)$, we then have $\bm{g}=[g_{(\betab)}]^\intercal$, $\bm{H}=[H_{(\bm{\beta},j)}]$ for $j=1,\cdots,M_i$ and $\betab=(\beta_{(1)},\cdots,\beta_{(d)})\in\mathbb{N}^d$, $0\leq|\betab|_{l_1}\leq 3$, where 
$$g_{(\betab)}=g_{(\beta_{(1)},\cdots,\beta_{(d)})}=\int_{B_\delta(\bm{0})}\gamma_\delta(|\xb-\yb|)(\yb-\xb_i)_{(1)}^{\beta_{(1)}}\cdots(\yb-\xb_i)_{(d)}^{\beta_{(d)}}d\yb \quad\text{ for }|\betab|_{l_1}\neq 0,$$
$$H_{(\bm{\beta},j)}=\gamma_{ji}\cdot(\xb_j-\xb_i)_{(1)}^{\beta_{(1)}}\cdots(\xb_j-\xb_i)_{(d)}^{\beta_{(d)}}.$$
By eliminating the constraints, the quadrature weights may be obtained by solving
\begin{equation}
\bm{\omega} = \bm{W}^{-1}\bm{H}^\intercal [\bm{H}\bm{W}^{-1}\bm{H}^\intercal]^{-1}\bm{g},
\end{equation}
where $\bm{H}\bm{W}^{-1}\bm{H}^\intercal=[m_{(\betab,\thetab)}]$ for $\betab=(\beta_{(1)},\cdots,\beta_{(d)})$, $\thetab=(\theta_{(1)},\cdots,\theta_{(d)})$, and $0\leq|\betab|_{l_1}\leq 3$, $0\leq|\thetab|_{l_1}\leq 3$ satisfying
$$m_{(\betab,\thetab)}=m_{(\thetab,\betab)},\quad m_{(\betab,\thetab)}=\sum_{j=1}^{M_i}\gamma_{ji}(\xb_j-\xb_i)^{\beta_{(1)}+\theta_{(1)}}_{(1)}\cdots(\xb_j-\xb_i)^{\beta_{(d)}+\theta_{(d)}}_{(d)}.$$
For problems where the reproducing constraints are redundant, $[\bm{H}\bm{W}^{-1}\bm{H}^\intercal]^{-1}$ may be replaced by the pseudo-inverse.

\subsubsection{Stability and Convergence Analysis}


To provide the stability proof of our method, we first show that the quadrature weights are all positive:
\begin{lem}\label{lem:weight}
 Consider a kernel $\gamma_\delta$ defined in \eqref{eqn:require_ga} with $s<d$.  For a given $\xb$ and sufficiently small $h/\delta$, quadrature weights obtained from \eqref{eq:quadQP} with the choice of $\bm{V}_h = \bm{S}_{\gamma_\delta,\xb}$ are all positive. In particular we have
\begin{equation}\label{eqn:wij}
\verti{\omega_{j,i}-h^d}\leq Ch^{\min(d+1,2d-s)},
\end{equation}
where $C$ is a constant independent of $h$.
\end{lem}

\begin{proof}
For a fixed $\delta$, we denote the weights as  $\omega^\delta_{j,i}$. Due to the scaling of the nonlocal kernel in \eqref{eqn:require_ga}, $\omega^\delta_{j,i}$ has the following scaling $\omega^\delta_{j,i}=\delta^d\omega^1_{j,i}$, and therefore it suffices to prove the estimate for $\delta=1$. With the symmetry property of $\gamma_1$ and the uniform Cartesian grid, we note that $g_{(\betab)}=0$ if one of the components of $\betab$ is an odd number, and $m_{(\betab,\thetab)}=0$ if one of the components of $\betab+\thetab$ is an odd number. Moreover, $g_{(\betab)}=g_{(\tilde{\betab})}$ if $\tilde{\betab}$ is a reordering of $\betab$, and $m_{(\betab,\thetab)}=m_{(\tilde{\betab},\tilde{\thetab})}$ if $\tilde{\betab}+\tilde{\thetab}$ is a reordering of ${\betab}+{\thetab}$. Denote $\bb=[\bm{H}\bm{W}^{-1}\bm{H}^\intercal]^{-1}\bm{g}$, we then have $b_{(\betab)}=0$ unless $|\betab|_{l_1}=0$ or $|\betab|_{l_1}=2$ {with only one nonzero entry $\beta_k=2$ for some $k\in\{1,\cdots,d\} $}. With the symmetry properties of $g_{(\betab)}$ and $m_{(\betab,\thetab)}$, we notice that $b_{(\betab)}=b_{(\tilde{\betab})}$ if $\tilde{\betab}$ is a reordering of $\betab$. Therefore, we may denote $b_{(0,\cdots,0)}:=b_0$, $b_{(2,0,\cdots,0)}=b_{(0,2,\cdots,0)}=\cdots=b_{(0,\cdots,0,2)}:=b_2$, and $b_0$, $b_2$ satisfy:
\begin{align*}
&\sum_{j=1}^{M_i}\gamma_{ji}b_0+\left[\sum_{j=1}^{M_i}\sum_{k=1}^d\gamma_{ji} (\xb_j-\xb_i)_{(k)}^2\right]\cdot b_2=\int_{B_1(\bm{0})}\gamma_1(|\yb|) d\yb,\\
&\left[\sum_{j=1}^{M_i}\gamma_{ji}(\xb_j-\xb_i)_{(1)}^2\right]b_0+\left[\sum_{j=1}^{M_i}\gamma_{ji}(\xb_j-\xb_i)_{(1)}^2\sum_{k=1}^d (\xb_j-\xb_i)_{(k)}^2\right]\cdot b_2=\int_{B_1(\bm{0})}\gamma_1(|\yb|)\yb_{(1)}^2 d\yb. 
\end{align*}
Notice that $s< d-1$ and compare the quadrature rule with the Riemann sum estimates for integrals, we obtain 
\begin{align*}
&\verti{\sum_{j=1}^{M_i}\gamma_{ji}h^d-\int_{B_1(\bm{0})}\gamma_1(|\yb|) d\yb}\leq C\int_{B_{\sqrt{d}h}(\bm{0})}\verti{\yb}^{-s} d\yb+Ch=C (h^{d-s}+h),\\
&\verti{h^d\sum_{j=1}^{M_i}\gamma_{ji}(\xb_j-\xb_i)_{(1)}^2-\int_{B_1(\bm{0})}\gamma_1(|\yb|)\yb_{(1)}^{2} d\yb}\leq C\int_{B_{\sqrt{d}h}(\bm{0})} \verti{\yb}^{2-s}d\yb+Ch=C(h^{d+2-s}+h)\leq Ch,\\
&\verti{h^d\sum_{j=1}^{M_i}\gamma_{ji}(\xb_j-\xb_i)_{(1)}^2\sum_{k=1}^d (\xb_j-\xb_i)_{(k)}^2-\int_{B_1(\bm{0})}\gamma_1(|\yb|)|\yb|^2\yb^2_{(1)} d\yb}\\
&~~~\leq C\int_{B_{\sqrt{d}h}(\bm{0})} \verti{\yb}^{4-s}d\yb+Ch= C(h^{d+4-s}+h)\leq Ch,
\end{align*}
where $B_{\sqrt{d}h}(\bm{0})$ is a sphere covering the origin point and the singularity on it. Note here we don't have an associated weight $\omega_{i,i}$, and the errors of integrals in the (hyper)cube containing the origin point with size $h^d$ were estimated separately. $B_{\sqrt{d}h}$ is then chosen as a sphere to contain this (hyper)cube. We obtain $b_0=h^d+O(h^{\min(d+1,2d-s)})$, $b_2=O(h^{\min(d+1,2d-s)})$, and $\bm{\omega}=\bm{W}^{-1}\bm{H}^\intercal\bb=h^d+O(h^{\min(d+1,2d-s)})$.
\end{proof}

From the above lemma and \cite{wendland2004scattered}, we can see that there exists a constant $C$ independent of $\delta$, such that for all $h<C\delta$ unisolvency holds for the optimization problem \eqref{eq:quadQP} and we have $\omega_{j,i}>0$. In the following we will denote this constant as $C_{pos}$ and generally require $h<C_{pos}{\delta}$. Note that the symmetry property in Lemma \ref{lem:weight} also yields $\omega_{j,i}=\omega_{i,j}$. 
Consider $A$ satisfying the conditions in \eqref{eqn:require_A} and 
\eqref{eqn:require_B} and $\gamma_\delta$ satisfying the conditions in Lemma \ref{lem:weight}, we obtain the discrete maximum principle for $\mcL_h^\delta$ and the well-poseness property of the discretized nonlocal diffusion problem \eqref{eq:discreteNonloc}:

\begin{lem}[Discrete Maximum/Minimum Principle]\label{lem:maxprinciple}
$v_h$ is a function on $\chi_h$ satisfying $\mcL_h^\delta v_h\geq0$ on $\chi_h\cap\omg$, then $\underset{\chi_h\cap\omg}{\max}v_h\leq \underset{\chi_h\cap\omg_\delta}{\max}v_h$. Similarly, when $v_h$ is a function on $\chi_h$ satisfying $\mcL_h^\delta v_h\leq0$ on $\chi_h\cap\omg$, then $\underset{\chi_h\cap\omg}{\min}v_h\geq \underset{\chi_h\cap\omg_\delta}{\min}v_h$.
\end{lem}
\begin{proof}
We first prove the maximum principle for $\mcL_h^\delta v_h\geq0$. Assume that $\underset{\chi_h\cap\omg}{\max}\,v_h\geq \underset{\chi_h\cap\omg_\delta}{\max}\,v_h$ and $v_h(\xb_i)=\underset{\chi_h\cap\omg}{\max}\,v_h$, then
\begin{align*}
0\geq&- \mcL_h^\delta v_h(\xb_i)=-2\sum_{\xb_j\in \chi_h\cap B_\delta(\xb_i)\backslash\{\xb_i\}}A(\xb_i,\xb_j)\gamma_\delta(|\xb_i-\xb_j|)\omega_{j,i}(u_h(\xb_j)-u_h(\xb_i))\geq 0.   
\end{align*}
Therefore, we have $u_h(\xb_j)=u_h(\xb_i)=\underset{\chi_h\cap\omg}{\max}u_h$. Apply the same argument to $\xb_j$ in $\omg$ and to their neighbors, one may obtain $u_h=const$ on $\chi_h$ and therefore  $\underset{\chi_h\cap\omg}{\max}u_h= \underset{\chi_h\cap\omg_\delta}{\max}u_h$. The discrete minimum principle can be similarly proved.
\end{proof}
With the discrete maximum principle,  the discretized nonlocal diffusion problem \eqref{eq:discreteNonloc} is therefore well-posed, i.e., there exists a unique solution to the discretized nonlocal diffusion problem \eqref{eq:discreteNonloc}.
We now consider the accuracy of the quadrature rule. In the following, we first present the truncation error estimate of the meshfree discretization for the nonlocal diffusion problem:
\begin{lem}\label{lem:consistency}
    Consider a kernel $\gamma_\delta$ satisfying the conditions in Lemma \ref{lem:weight}, $A(\yb,\xb)$, $u(\xb)$ are $C^{1}$ with respect to $\xb$, and a fixed $\delta$. Then for $h<C_{pos}\delta$ quadrature weights obtained from \eqref{eq:quadQP} with the choice of $\bm{V}_h =  \bm{S}_{\gamma_\delta,\xb}$ satisfy the following pointwise error estimate, with $C>0$ independent of $h$:
{\begin{equation}\label{eqn:consistency}
\underset{\xb_i\in\chi_h\cap\omg}{\max}\verti{\mcL^\delta[u](\xb_i)-\mcL_h^\delta[u](\xb_i)}\leq Ch^{\min(1,d-s)}.    
\end{equation}}
\end{lem}

\begin{proof}
We note that
\begin{align*}
&\sum_{j=1}^{M_i} \verti{A(\xb_i,\xb_j)\gamma_\delta(|\xb_i-\xb_j|) (u(\xb_j)-u(\xb_i))}h^{\min(d+1,2d-s)}\leq C\sum_{j=1}^{M_i} \verti{\xb_i-\xb_j}^{1-s}h^{\min(d+1,2d-s)}\\
&\leq C\int_{B_\delta(\bm{0})} \verti{\yb}^{1-s}d\yb\,h^{\min(1,d-s)}\leq C(\delta)h^{\min(1,d-s)}.
\end{align*}
On the other hand, {let $G(\xb_j)$ denote the (hyper)cubic of size $h^d$ centered at $\xb_j$ }, we have the error estimate of the Riemann sum formulation as:
\begin{align*}
&\left|\int_{B_\delta(\xb_i)} A(\xb_i,\yb)\gamma_\delta(|\xb_i-\yb|) (u(\yb)-u(\xb_i))\,d\yb - \sum_{j=1}^{M_i} A(\xb_i,\xb_j)\gamma_\delta(|\xb_i-\xb_j|) (u(\xb_j)-u(\xb_i)) h^d\right|\\
&\leq Ch+C\int_{B_{\sqrt{d}h}(\bm{0})}|\yb|^{1-s}d\yb\\
&+\sum_{j=1}^{M_i} \left|\int_{G(\xb_j)} A(\xb_i,\yb)\gamma_\delta(|\xb_i-\yb|) (u(\yb)-u(\xb_i))\,d\yb - A(\xb_i,\xb_j)\gamma_\delta(|\xb_i-\xb_j|) (u(\xb_j)-u(\xb_i)) h^d\right|\\
&\leq C(h+h^{d-s+1})+C \sum_{j=1}^{M_i} h^{d+1} \max_{\zb\in G(\xb_j)}|\nabla_{\yb} (A(\xb_i,\zb)(u(\zb)-u(\xb_i))\gamma_\delta(|\xb_i-\zb|))|\\
&\leq Ch+C \sum_{j=1}^{M_i} h^{d+1} \left(\max_{\zb\in G(\xb_j)}|(u(\zb)-u(\xb_i))\nabla_{\yb}\gamma_\delta(|\xb_i-\zb|)|+\max_{\zb\in G(\xb_j)}\gamma_\delta(|\xb_i-\zb|)|\nabla_{\yb} (A(\xb_i,\zb)(u(\zb)-u(\xb_i))|\right)\\
&\leq Ch+C \sum_{j=1}^{M_i} h^{d+1} \left(\max_{\zb\in G(\xb_j)}|\xb_i-\zb|^{-s}\right)\leq Ch+Ch\int_{B_{\delta}(\bm{0})}|\yb|^{-s}d\yb\leq Ch
\end{align*}
where $\verti{\nabla_{\yb} a(\xb,\yb)}$ denotes the maximum component of the first order partial derivatives with respect to $\yb$. Therefore, substituting \eqref{eqn:wij} into \eqref{eq:discreteNonloc} yields
\begin{align*}
&\verti{\mcL^\delta[u](\xb_i)-\mcL_h^\delta[u](\xb_i)}\\
=&\left|\int_{B_\delta(\xb_i)} A(\xb_i,\yb)\gamma_\delta(|\xb_i-\yb|) u(\yb)\,d\yb - \sum_{j=1}^{M_i} A(\xb_i,\xb_j)\gamma_\delta(|\xb_i-\xb_j|) u(\xb_j) h^d\right|+O(h^{\min(1,d-s)})\\
\leq&Ch^{\min(1,d-s)}.
\end{align*}
Note here the constant $C$ is independent of $h$ but may depends on $\delta$.
\end{proof}

To prove the asymptotic compatibility, we need the truncation error estimate to the local limit. In \cite[Theorem~2.1]{trask2019asymptotically}, the authors have shown that for a sufficiently smooth $u$, the quadrature weights obtained from \eqref{eq:quadQP} with the choice of $\bm{V}_h =  \bm{S}_{\gamma_\delta,\xb}$ provides an $O(\delta^{2})$ pointwise error bound for the integral approximation of $2\int_{ B_\delta(\xb_i)}\gamma_\delta(|\xb_i-\yb|) (u(\yb)-u(\xb_i))\,d\yb$ when the ratio $h/\delta$ is fixed. We can easily extend this error estimate to a nonlocal diffusion operator with heterogeneous diffusion coefficient:

\begin{lem}\label{lem:AC_1}
  Consider a kernel $\gamma_\delta$ satisfying the conditions in Lemma \ref{lem:weight}, $A(\xb,\yb)$ and $u(\yb)$ are $C^{4}$ with respect to $\yb$, and fixed ratio $h/\delta<C_{pos}$. Quadrature weights obtained from \eqref{eq:quadQP} with the choice of $\bm{V}_h = \bm{S}_{\gamma_\delta,\xb}$ satisfy the following pointwise error estimate, with $C>0$ independent of $\delta$:
$$\underset{\xb_i\in\bm{V}_h}{\max}\verti{\mcL^\delta[u](\xb_i)-\mcL_h^\delta[u](\xb_i)}\leq C \delta^2.$$
\end{lem}
\begin{proof}
Note that for a given $\xb$ we consider the estimate of $A(\xb,\yb)u(\xb)\in C^4$ and  $A(\xb,\yb)u(\yb)\in C^4$, the following truncation estimate is obtained immediately from \cite[Theorem~2.1]{trask2019asymptotically}:
\begin{align*}
&\underset{\xb_i\in\chi_h\cap\omg}{\max}\left|2\int_{ B_\delta(\xb_i)}\frac{A(\xb_i,\yb)D_0}{|\xb_i-\yb|^{s}} (u(\yb)-u(\xb_i))\,d\yb\right.\\
&~~~~~~\left.-2\sum_{\xb_j\in\chi_h\cap B_\delta(\xb_i)\backslash\{\xb_i\}}\omega_{j,i}\frac{A(\xb_i,\xb_j)D_0}{|\xb_i-\xb_j|^{s}} (u(\xb_j)-u(\xb_i))\right|\leq C\delta^{4-s+d}.
\end{align*}
And the proof is finished by taking $\gamma_\delta(|\xb-\yb|)=\frac{D_0}{\delta^{d+2-s}|\xb_i-\yb|^{s}}$ as stated in \eqref{eqn:require_ga}.
\end{proof}

From the derivation of Theorem \ref{thm:compatibility}, one can see that $\verti{\mcL^0[u](\xb_i)-\mcL^\delta[u](\xb_i)}=O(\delta^2)$. We therefore obtain the following truncation error estimate to the local limit:
\begin{cor}\label{lem:AC}
  Consider a kernel $\gamma_\delta$ satisfying the conditions in Lemma \ref{lem:weight}, $A(\xb,\cdot),u(\cdot) \in C^4$ and fixed ratio $h/\delta<C_{pos}$. Quadrature weights obtained from \eqref{eq:quadQP} with the choice of $\bm{V}_h =  \bm{S}_{\gamma_\delta,\xb}$ satisfy the following pointwise error estimate, with $C>0$ independent of $\delta$:
$$\underset{\xb_i\in\bm{V}_h}{\max}\verti{\mcL^0[u](\xb_i)-\mcL_h^\delta[u](\xb_i)}\leq C \delta^2.$$
\end{cor}

With the discrete maximum principle and the above truncation estimate results, we finally get the main results on the stability in maximum norm:
\begin{lem}[Stability]\label{lem:discstable}
Consider a bounded domain $\omg$ 
and $\gamma_\delta$ satisfying the conditions in Lemma \ref{lem:weight}. Assume that $a\in C^{\infty}(\omg)$ and $\omg\cup\omg_\delta\in C^1$, then there exist generic constants $C$ and $\delta_0>0$ such that when $\delta<\delta_0$ and $h<C_{pos}\delta$, solution to the discretized nonlocal diffusion problem \eqref{eq:discreteNonloc} satisfies:
$$\underset{\xb_i\in\chi_h}{\max}|u^\delta_h(\xb_i)|\leq \underset{\omg_\delta}{\max}\,|u^D|+C\underset{\xb_i\in\chi_h\cap \omg}{\max}|f(\xb_i)|.$$
Here $C$ is independent of both $h$ and $\delta$.
\end{lem}

\begin{proof}
We first construct a barrier function $\psi(\xb)$. With the properties of $a$, there exists a $C^\infty$ solution $\hat{\psi}$ for the following classical diffusion problem \cite{Evans02}:
\begin{align*}
\mcL^0 \hat{\psi}=2 & \quad \text{ on }\omg\cup\omg_\delta,\\
\hat{\psi}=0& \quad \text{ on }\partial(\omg\cup\omg_\delta),    
\end{align*}
where $\partial(\omg\cup\omg_\delta)$ denotes the exterior boundary of $\omg\cup\omg_\delta$. We then define $\psi(\xb):=\hat{\psi}(\xb)-\underset{\zb\in\omg\cup\omg_\delta}{\min}\hat{\psi}(\zb)$, and notice that $\psi$ satisfies: $\mcL^0\psi=2>0$, $0\leq\psi\leq \left(\underset{\zb\in\omg\cup\omg_\delta}{\max}\hat{\psi}(\zb)-\underset{\zb\in\omg\cup\omg_\delta}{\min}\hat{\psi}(\zb)\right):=M_\psi$. Set $M_f:=\underset{\xb_i\in\chi_h\cap \omg}{\max}|f(\xb_i)|$. Since $\underset{\xb_i\in\bm{V}_h}{\max}\verti{\mcL^0[{\psi}](\xb_i)-\mcL_h^\delta[{\psi}](\xb_i)}\leq C \delta^2$ as shown in Corollary \ref{lem:AC}, there exists a constant $\delta_0>0$ such that when $\delta<\delta_0$ we have $\mcL_h^\delta\psi\geq1>0$. Then $\mcL_h^\delta(u^\delta_h+M_f\psi)(\xb_i)=-f(\xb_i)+M_f\geq 0$, and the discrete maximum principle yields
\begin{align}
\nonumber&\underset{\xb_i\in\chi_h}{\max}u^\delta_h(\xb_i)\leq \underset{\xb_i\in\chi_h}{\max}(u^\delta_h+M_f\psi)(\xb_i)\leq \underset{\xb_i\in\chi_h\cap\omg_\delta}{\max}(u^\delta_h+M_f\psi)(\xb_i)\\
&\leq \underset{\omg_\delta}{\max}\,|u^D|+M_fM_\psi=\underset{\omg_\delta}{\max}\,|u^D|+M_\psi\underset{\xb_i\in\chi_h\cap \omg}{\max}|f(\xb_i)|.\label{eqn:maxuhd}
\end{align}
Similarly we can show that $\mcL_h^\delta(u^\delta_h-M_f\psi)(\xb_i)=-f(\xb_i)-M_f\leq 0$ and 
$$\underset{\xb_i\in\chi_h}{\min}u_h(\xb_i)\geq -\underset{\omg_\delta}{\max}\,|u^D|-M_\psi\underset{\xb_i\in\chi_h\cap \omg}{\max}|f(\xb_i)|,$$
which together with \eqref{eqn:maxuhd} finishes the proof.
\end{proof}


With the stability property of $\mcL_h^\delta$ and the truncation estimates in Lemmas \ref{lem:consistency}-\ref{lem:AC} we proceed to prove the two types of convergence results in \eqref{eqn:convh}. We first consider the case with fixed $\delta$ and vanishing $h$. In particular, we investigate the convergence of numerical solution to nonlocal solution by combining Lemma \ref{lem:consistency} with the stability property:
\begin{thm}[Convergence to Deterministic Nonlocal Solution]\label{thm:consistency}
Assume that the conditions in Lemma \ref{lem:consistency} and Lemma \ref{lem:discstable} are satisfied and $u^\delta\in C^1(\omg)$, then there exists a $\delta_0>0$ such that for a fixed $\delta$ satisfying $0<\delta<\delta_0$ and $h<C_{pos}\delta$, the following convergence property holds for the numerical solution of \eqref{eq:discreteNonloc}:
\begin{equation}
\vertii{u_h^\delta-u^\delta}_{L^\infty(\chi_h)}\leq C h^{\min(1,d-s)},
\end{equation}
where $C$ is a generic constant independent of $h$.
\end{thm}
\begin{proof}
Apply the stability theorem to $u_h^\delta-u^\delta$, we immediately obtain
$$\underset{\xb_i\in\chi_h}{\max}|u^\delta_h(\xb_i)-u^\delta(\xb_i)|\leq C\underset{\xb_i\in\chi_h\cap \omg}{\max}|\mcL^\delta_h u^\delta(\xb_i)-\mcL^\delta u^\delta(\xb_i)|\leq C h^{\min(1,d-s)}.$$
\end{proof}

Next, we show the AC property of the meshfree method, when both $\delta$ and $h$ vanish with a fixed ratio $\delta/h$. In particular, we investigate the convergence rate of numerical solution to local limit by combining Lemma \ref{lem:AC_1} and Corollary \ref{lem:AC} with the stability property of $\mcL_h^\delta$:
\begin{thm}[Asymptotic Compatibility in Deterministic Nonlocal Problems]\label{thm:AC}
Assume that the conditions in Corollary \ref{lem:AC} and Lemma \ref{lem:discstable} are satisfied, and $u^\delta,u^0\in C^4(\omg)$, then there exists a $\delta_0>0$ such that for any $0<\delta<\delta_0$ and fixed ratio $h/\delta<C_{pos}$, the following error estimate holds for $u^\delta_h$:
\begin{equation}
\vertii{u_h^\delta-u^\delta}_{L^\infty(\chi_h)}\leq C\delta^2.
\end{equation}
 Moreover, the meshfree scheme \eqref{eq:discreteNonloc} is asymptotically compatible, i.e.,
\begin{equation}
\vertii{u_h^\delta-u^0}_{L^\infty(\chi_h)}\leq C\delta^2,
\end{equation}
where $C$ is a generic constant independent of $\delta$.
\end{thm}
\begin{proof}
Apply the stability theorem to $u_h^\delta-u^\delta$ and $u_h^\delta-u^0$, we immediately obtain:
$$\underset{\xb_i\in\chi_h}{\max}|u^\delta_h(\xb_i)-u^\delta(\xb_i)|\leq C\underset{\xb_i\in\chi_h\cap \omg}{\max}|\mcL^\delta_h u^\delta(\xb_i)-\mcL^\delta u^\delta(\xb_i)|\leq C\delta^2,$$
$$\underset{\xb_i\in\chi_h}{\max}|u^\delta_h(\xb_i)-u^0(\xb_i)|\leq C\underset{\xb_i\in\chi_h\cap \omg}{\max}|\mcL^\delta_h u^0(\xb_i)-\mcL^0 u^0(\xb_i)|\leq C\delta^2.$$
\end{proof}

\begin{rem}
When taking a flat kernel with $s=0$, from Theorems \ref{thm:consistency} and \ref{thm:AC}, we note that the optimal convergence rate to the nonlocal solution is $O(h)$ when $h\rightarrow0$, while the optimal convergence to the local limit is $O(\delta^2)$ when $h,\delta\rightarrow0$.
\end{rem}

\subsection{Stochastic: Probabilistic Collocation Method with Sparse Grids}

 To solve the stochastic problem introduced in Section \ref{sec:randomdiff}, we employ the probabilistic collocation method (PCM) in the parametric space for its high resolution and ease of implementation by sampling at discrete points in random space \cite{tatang1994direct,keese2003numerical,xiu2005high}. Consider the stochastic equation \eqref{eqn:nonlocal_random}, PCM can be seen as a Lagrange interpolation in the random space. In particular, let $\Theta_N=\{\xib_k\}_{k=1}^K\subset \Gamma$ be a set of prescribed nodes such that the Lagrange interpolation in the random space $\Gamma$ is poised in an interpolation space $\Gamma_I$, where $N$ is the dimension of the parametric space. Then any function $v:\Gamma\rightarrow \real$ can be approximated using the Lagrange interpolation polynomial:
$$\mathcal{J}[v](\xib)=\sum_{k=1}^K v(\xib_k)J_k(\xib), $$
where $J_k(\xib)$ is the Lagrange polynomial satisfying $J_k(\xib)\in \Gamma_I$ and $J_k(\xib_j)=\delta_{kj}$. Denoting $\hat{u}(\xb,\xib):=\sum_{k=1}^K u(\xb,\xib_k)J_k(\xib)$, the collocation procedure to solve the stochastic nonlocal equation is
$$R(\hat{u}(\xb,\xib))|_{\xib_k}=0,\quad\forall k=1,\cdots,K,$$
where $R$ is the residual of \eqref{eqn:nonlocal_random}. With the property of Lagrange interpolation, we obtain
\begin{equation}\label{eqn:lag}
\left\{\begin{array}{ll}
-\mcL^\delta[u](\xb,\xib_k)=-\mathcal{D}[A(\xb,\yb,\xib_k)\mathcal{G}[u(\xb,\xib_k)]]=f(\xb)&\text{ for }\xb\in\omg,\\
u(\xb,\xib_k)=u^D(\xb),\quad &\text{ for }\xb\in\omg_\delta,
\end{array}\right.
\end{equation}
for $k=1,\cdots,K$. Note that \eqref{eqn:lag} is equivalent to solving $K$ deterministic nonlocal diffusion problems, where the deterministic meshfree solver discussed in Section \ref{sec:meshfree} can be readily applied. Therefore, the PCM approach can be implemented in an embarrassingly parallel way and the total computational cost is the product of the number of collocation points times the cost of the deterministic problem.

With the numerical solution of \eqref{eqn:lag} on all collocation points $\xib_i$, the statistical moments of the random solution can be evaluated:
$$\mathbb{E}[u](\xb)\approx \mathbb{E}[\hat{u}](\xb)=\int_\Gamma \sum_{k=1}^K u(\xb,\xib_k)J_k(\xib)\rho(\xib)d\xib,$$
$$\sigma[u](\xb)\approx \sigma[\hat{u}](\xb)=\sqrt{\int_\Gamma \left[\sum_{k=1}^K u(\xb,\xib_k)J_k(\xib)\right]^2\rho(\xib)d\xib-[\mathbb{E}[\hat{u}](\xb)]^2},$$
and so on. Here $\rho$ is the PDF of random variable $\xib$. To further approximate the integral for above polynomials, we employ the quadrature rule approximation by choosing the set $\Theta_N$ as quadrature point set:
\begin{align}
&\mathbb{E}[u](\xb)\approx \mathbb{E}[\hat{u}](\xb)\approx \sum_{k=1}^K u(\xb,\xib_k)\mu_k,\label{eqn:E}\\
&\sigma[u](\xb)\approx \sigma[\hat{u}](\xb)\approx\sqrt{\sum_{k=1}^K u^2(\xb,\xib_k)\mu_k-[\mathbb{E}(\hat{u})(\xb)]^2},\label{eqn:std}
\end{align}
where $\{\mu_k\}_{k=1}^K$ is the set of corresponding quadrature weights.

There are mainly two different strategies for the selection of collocation point sets: the tensor products of 1D collocation point sets and a sparse grid strategy for high dimensionality. In the tensor product strategy, one first construct a 1D interpolation for each dimension in the random space. For the $i$-th dimension, we take $\varpi_{(i)}$ numbers of nodal points $\Theta^{\varpi_{(i)}}_1=\{\xi_1^i,\cdots,\xi_{\varpi_{(i)}}^i\}\subset[-1,1]$, a 1D interpolation for a smooth function $v$ on the $i$-th dimension then writes:
\begin{equation}
  \mathcal{U}^{\varpi_{(i)}}[v](\xi_{(i)})=\sum_{k=1}^{\varpi_{(i)}}v(\xi^i_k)J^i_k(\xi_{(i)})
\end{equation}
where $J^i_k(\xi_{(i)})$ is the 1D Lagrange polynomial. Then for the multivariate case $v:\real^N\rightarrow\real$, the tensor product formula is:
\begin{equation}\label{eqn:tensor}
    \mathcal{J}[v]=\left(\mathcal{U}^{\varpi_{(1)}}\otimes\cdots\otimes\mathcal{U}^{\varpi_{(N)}}\right)[v]=\sum_{k_1=1}^{\varpi_{(1)}}\cdots\sum_{k_N=1}^{\varpi_{(N)}} v\left(\xi^1_{k_1},\cdots,\xi^N_{k_N}\right)\left(J^1_{k_1}\otimes\cdots\otimes J^N_{k_N}\right).
\end{equation}
Notice here \eqref{eqn:tensor} requires $K=\Pi_{i=1}^{N}\varpi_{(i)}$ numbers of collocation points in total, which grows quickly when $N$ gets large. Therefore, the tensor product strategy may be employed for problems with a small number of random dimensions, but its required number of collocation points $K$ generally grows exponentially as $N$ increases and makes the simulation non-feasible (see, e.g., \cite{lin2009efficient}). Hence for problems with a relatively large random dimension, we employ the sparse grid strategy. In particular, we employ the sparse grid constructed by the Smolyak algorithm \cite{smolyak1963quadrature}, which is a linear combination of tensor product formulas:
\begin{equation}\label{eqn:sparsegrid}
  \mathcal{J}[v]=\sum_{\zeta-N+1\leq|\bm{\varpi}|\leq \zeta}(-1)^{\zeta-\verti{\bm{\varpi}}_{l_1}}\binom{N-1}{\zeta-\verti{\bm{\varpi}}_{l_1}}\left(\mathcal{U}^{\varpi_{(1)}}\otimes\cdots\otimes\mathcal{U}^{\varpi_{(N)}}\right).
\end{equation}
Here $\zeta$ is the sparseness parameter, $\bm{\varpi}=(\varpi_{(1)},\cdots,\varpi_{(N)})\in\mathbb{N}^N$, $\verti{\bm{\varpi}}_{l_1}=\sum_{i=1}^N \varpi_{(i)}$, and $\varpi_{(i)}$ represents the number of collocation points in random dimension $i$. To compute \eqref{eqn:sparsegrid}, only evaluations on the sparse grids are needed:
\begin{equation}\label{eqn:sparsesample}
\Theta_N=\underset{\zeta-N+1\leq|\bm{\varpi}|_{l_1}\leq \zeta}{\bigcup} \left(\Theta^{\varpi_{(1)}}_1\times\cdots\times\Theta^{\varpi_{(N)}}_1\right).
\end{equation}
As shown in \cite{novak1996high,novak1999simple}, \eqref{eqn:sparsegrid} is exact for $p(\xib)\in\mathbb{P}_{\zeta-N}(\real^N)$ (all polynomials of degree less than $\zeta-N$) and the total number of nodes $K\sim \frac{{2N}^{\zeta-N}}{(\zeta-N)!}$. Therefore, we may see that the sparse grid formulation typically requires a much smaller number of collocation points $K$ than the full tensor product set and we will refer $\eta=\zeta-N$ as the ``level'' of the Smolyak formulation.

\subsection{An Asymptotically Compatible Meshfree PCM}

Let $\delta$ be the horizon size, $h$ be the grid size in space, and $K=K(\eta, N)$ be the total number of collocation points we use in the parameter space. We then denote the numerical solution to \eqref{eqn:nonlocal_random} by $u^\delta_{h,K}$. As $h\to0$ and $K\to \infty$, we expect the numerical solution to converge to the exact solution $u^\delta$ for \eqref{eqn:nonlocal_random} in $L^2(\omg)\otimes L^2_\rho(\Gamma)$. In addition, our numerical method is asymptotically compatible, i.e., $u^\delta_{h,K}\to u^0$ as $h\to0$, $K\to \infty$, and $\delta\to0$, where $ u^0$ is the exact solution for \eqref{eqn:local_random}. 
The error estimate can be facilitated by introducing an intermediate function 
$u^\delta_K$, the semi-discrete solution to \eqref{eqn:nonlocal_random}. Then we split the errors into
\[
\left\{
\begin{aligned}
&u^\delta - u^\delta_{h,K} = (u^\delta - u^\delta_K) + (u^\delta_K -u^\delta_{h,K}) \\
&u^0 - u^\delta_{h,K} = (u^0 - u^\delta)+(u^\delta-u^\delta_K) + (u^\delta_K -u^\delta_{h,K}). 
\end{aligned}
\right.
\]
For each $\xib\in \Gamma$, 
\[
u^\delta_K(\xib) =\sum_{k=1}^K u^\delta(\xib_k)J_k(\xib),
\text{ and }
u^\delta_{h,K}(\xib) =\sum_{k=1}^K u^\delta_h(\xib_k)J_k(\xib). 
\]
The estimates of $\| u^\delta(\xib_k) - u^\delta_{h}(\xib_k)\|_{L^2(\omg)}$ for each $\xib_k$ 
are followed by the estimates in Section \ref{sec:method}. So we have 
$\| u^\delta_K(\xib) - u^\delta_{h,K}(\xib) \|_{L^2(\omg)}\to 0$ for each $\xib\in\Gamma$ 
and therefore $\|u^\delta_K- u^\delta_{h,K}  \|_{L^2(\omg)\otimes L^2_\rho(\Gamma)}\to0$ as $h\to0$. The rate at which $\|u^\delta_K- u^\delta_{h,K}  \|_{L^2(\omg)\otimes L^2_\rho(\Gamma)}$ converges to zero follows the estimates in Section \ref{sec:method}.  Moreover, we have 
$\| u^0(\xib) - u^\delta(\xib)\|_{L^2(\omg)}\leq C \delta^2$ for any $\xib$ as long as $u^0\in C^4(\omg)$. 
Therefore we only need to estimate $u^\delta - u^\delta_{K}$.

From the analytic regularity, $\xib\mapsto u^\delta(\xib)$ is a $S_\delta(\omg)$-valued map that admits an analytic extension to an open neighborhood $\mcA\subset \C^N$ of $\Gamma$. Moreover, we know that $\max_{\hat\xib \in \mcA} \| u^\delta(\hat\xib)\|_{S_\delta(\omg)}$ is uniformly bounded and independent of $\delta$. Therefore, the estimate of $u^\delta - u^\delta_{K}$ is followed exactly from \cite{nobile2008anisotropic}. Here we present the following result, whose proof can be found in \cite[Theorems 3.10-3.11]{nobile2008anisotropic}. 
\begin{lem}
\label{lem:convergenceinK}
 There exists $C_1>0$ and $\beta_1>0$ depending on $N$ and the analytic region $\mcA$ such that  
 \begin{equation}
 \label{eqn:algebraic}
\max_{\xib\in\Gamma}\|u^\delta(\xib) - u^\delta_K(\xib)\|_{S_\delta(\omg)}  \leq C_1 K^{-\beta_1}.
 \end{equation}
 Moreover, when $\eta > \frac{N}{\log(2)}$, there exists $C_2>0$, $C_3>0$ and $\beta_2>0$ depending on $N$ and the analytic region $\mcA$, and $\beta_3>0$ depending only on $N$ such that 
\begin{equation}
 \label{eqn:subexponential}
\max_{\xib\in\Gamma}\|u^\delta(\xib) - u^\delta_K(\xib)\|_{S_\delta(\omg)}  \leq C_2 K^{\beta_2} e^{- C_3 K^{\beta_3}}.
 \end{equation}
\end{lem}
Equation \eqref{eqn:algebraic} shows at least algebraic convergence with respect to the number of collocation points $K$, while equation \eqref{eqn:subexponential} shows the subexponential convergence when the level of Smolyak formulation $\eta >\frac{N}{\log(2)} $. In the numerical experiments below, we will choose large enough $\eta$ so as to observe the subexponential convergence. Based on Lemma \ref{lem:convergenceinK} and the above discussions, we have the following convergence theorems. 

\begin{thm}\label{thm:nonlocalconv}
Assume that $\max_{\xib\in\Gamma} \| u^\delta(\xib) \|_{C^p(\omg)} < \infty$ where $p\geq 1+s$,  then there exists a $\delta_0>0$ such that for a fixed $\delta$ satisfying $0<\delta<\delta_0$ and any  $h<C_{pos}\delta$, we have
\[
\| u^\delta - u^\delta_{h, K}\|_{L^2(\omg)\otimes L^2_\rho(\Gamma)} \leq C h^{\min(1,d-s)}+ C_1 K^{-\beta_1}.
\]
Moreover, if
$\eta >\frac{N}{\log(2)} $, we have
\[
\| u^\delta - u^\delta_{h, K}\|_{L^2(\omg)\otimes L^2_\rho(\Gamma)} \leq C h^{\min(1,d-s)} + C_2 K^{\beta_2} e^{- C_3 K^{\beta_3}}.
\]
The constants $C_1$, $C_2$, $C_3$, $\beta_1$, $\beta_2$, $\beta_3$ are defined in Lemma \ref{lem:convergenceinK}. 
\end{thm}
\begin{proof}
Let $(u^\delta - u^\delta_K) + (u^\delta_K -u^\delta_{h,K})$, then the theorem is a combination of Theorem \ref{thm:consistency} and Lemma \ref{lem:convergenceinK}. 
\end{proof}

\begin{thm}\label{thm:localconv}
Let $\delta_0$ be the constant defined in Lemma \ref{lem:discstable}. 
Assume that $\max_{\delta\in [0,\delta_0)}\max_{\xib\in\Gamma} \| u^\delta(\xib) \|_{C^4(\omg)} < \infty$. Then for any $\delta$ with  $0<\delta<\delta_0$ and fixed ratio $h/\delta<C_{pos}$, we have 
\[
\| u^0 - u^\delta_{h, K}\|_{L^2(\omg)\otimes L^2_\rho(\Gamma)} \leq C \delta^2+ C_1 K^{-\beta_1}.
\]
Moreover, if
$\eta >\frac{N}{\log(2)} $, we have
\[
\| u^0 - u^\delta_{h, K}\|_{L^2(\omg)\otimes L^2_\rho(\Gamma)} \leq C \delta^2+ C_2 K^{\beta_2} e^{- C_3 K^{\beta_3}}.
\]
The constants $C_1$, $C_2$, $C_3$, $\beta_1$, $\beta_2$, $\beta_3$ are defined in Lemma \ref{lem:convergenceinK}.
\end{thm}
\begin{proof}
Let $u^0 - u^\delta_{h,K} = (u^0 - u^\delta)+(u^\delta-u^\delta_K) + (u^\delta_K -u^\delta_{h,K})$. 
First we know that $\|u^0-u^\delta \|\leq C\delta^2$. Moreover, from Theorem \ref{thm:AC} we have $\| u^\delta_K -u^\delta_{h,K}\|\leq C \delta^2 $.  Then the theorem is combination of these estimates and Lemma \ref{lem:convergenceinK}.
\end{proof}



\section{Numerical Verification of Convergences}\label{sec:veri}

In this section, we numerically verify the proposed approach by investigating the two types of convergences: the consistency to nonlocal solutions and then the asymptotic compatibility to local companions\footnote{Notice that although we assumed $\Gamma= \prod_{i=1}^N [-1,1]\subset \R^N$ in Sections \ref{sec:pde}-\ref{sec:method}, in numerical tests we investigate and numerically verify our analysis on more general cases.}. In particular, we study the $L^2$ errors for the first two statistical moments, the mean and standard deviation (std). Let $u^{\delta}_{h,K}$ represent the numerical solution with spatial grid size $h$ in meshfree methods and $K$ samples in sparse grid PCM, $u^\delta$ represents the analytical nonlocal solution and $u^0$ stands for the analytical local limit, in Section \ref{sec:nonlocaltest} we investigate the convergence of numerical solutions to
the nonlocal analytical solution with vanishing $h$ and increasing sample numbers $K$ by calculating:
\begin{equation}
\| {\mathbb{E}}(u^{\delta}_{h,K})- \mathbb{E}(u^\delta) \| _{L_2(\omg)}, \quad \text{ and }\quad\|{\sigma}(u^{\delta}_{h,K})- \sigma(u^\delta) \|_{L_2(\omg)},
\end{equation}
where the mean and standard deviation of $u^{\delta}_{h,K}$ are numerically evaluated with the quadrature rule approximation in \eqref{eqn:E}-\eqref{eqn:std}. Similarly, in Section \ref{sec:localtest} we investigate the convergence of numerical solutions to
the local limit as $K$ increases and $\delta,h$ goes to $0$ simultaneously under the $\delta$-convergence limit. In particular we calculate:
\begin{equation}
\| {\mathbb{E}}(u^{\delta}_{h,K})- \mathbb{E}(u^0) \| _{L_2(\omg)}, \quad \text{ and }\quad\|{\sigma}(u^{\delta}_{h,K})- \sigma(u^0) \|_{L_2(\omg)}.
\end{equation}



\subsection{Consistency with the Stochastic Nonlocal Solution}\label{sec:nonlocaltest}

In this section we study the consistency of the numerical solution to the nonlocal analytical solution, on both 1D and 2D physical domains.

\subsubsection*{Test 1: consistency study on a problem with 1D physical domain and 5D parametric space}

\begin{figure}[h!]
\centering
\subfigure[Convergence with refinement in the physical space.]{\includegraphics[width=.45\columnwidth]{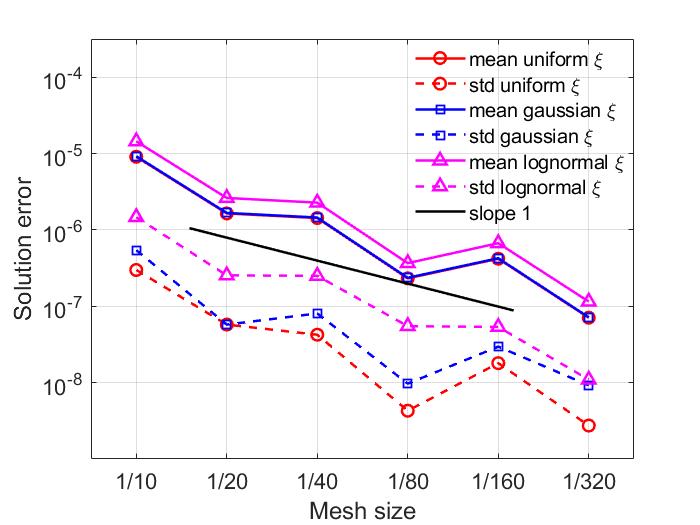}}\\
\subfigure[Convergence with sample numbers in the log scale.]{\includegraphics[width=.45\columnwidth]{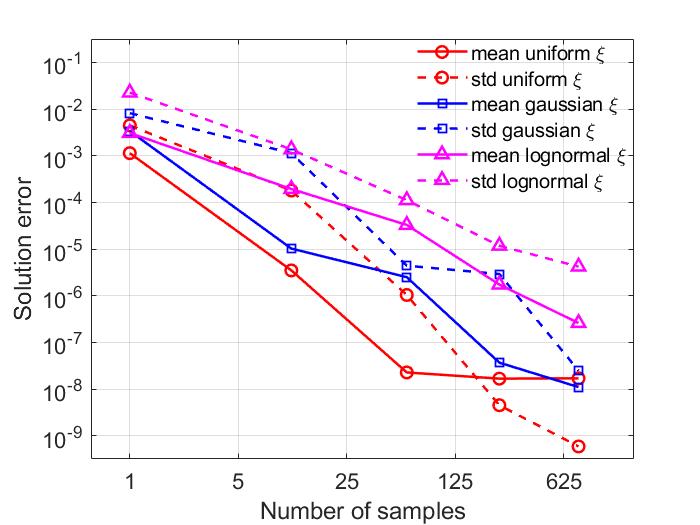}}
\subfigure[Convergence with sample numbers in the linear scale.]{\includegraphics[width=.45\columnwidth]{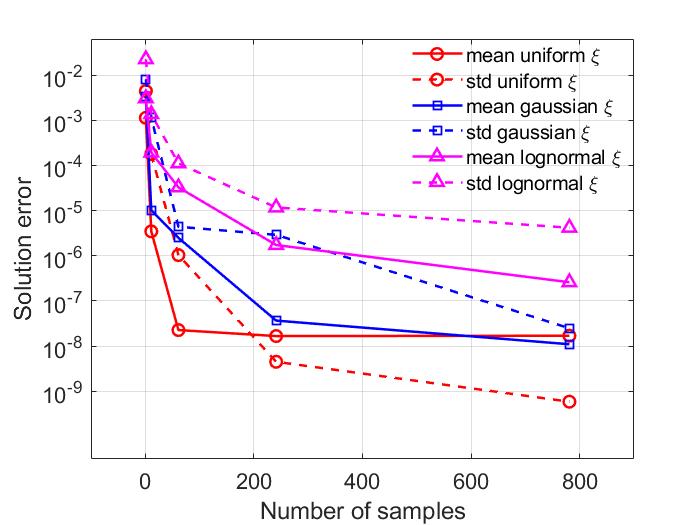}}
\caption{Test 1: consistency study of the numerical solution to the nonlocal analytical solution on a problem with 1D physical domain and 5D parametric space. Results in (a) are generated with $781$ samples, which corresponds to Smolyak formulation level $5$. The data points in (b) and (c) correspond to Smolyak formulation levels $\eta=1,\cdots,5$. }
\label{fig:test11,ex2it4}
\end{figure}

We consider a case with 1D physical domain $\Omega=[-1,1]$ and 5D parametric space $\xib=(\xi_{(1)},\cdots,\xi_{(5)})$, where $\xi_{(i)}$ are i.i.d. random variables. We manufacture the nonlocal analytical solution 
$$u^\delta(x,\xib)=\cos(0.5x)/(5+\cos(\xi_{(1)})+\sin(2\xi_{(2)})+\cos(3\xi_{(3)})+\sin(4\xi_{(4)})+\cos(5\xi_{(5)})),$$ 
with fixed $\delta = 0.38$, nonlocal diffusion coefficient 
$$A(x,y,\xib)=(2+\cos(0.5(x+y)))(5+\cos(\xi_{(1)})+\sin(2\xi_{(2)})+\cos(3\xi_{(3)})+\sin(4\xi_{(4)})+\cos(5\xi_{(5)})),$$
and loading
\begin{equation}\nonumber
\begin{split}
f(x) = -\frac{3}{\delta^3}\Big( \frac{1}{2}\big(-2\sin(0.5(3x+\delta))+\sin(0.5(3x+2\delta))-3\delta \cos(0.5x)+6\sin(0.5(x+\delta))\big) \\ -\frac{1}{2}\big(-2\sin(0.5(3x-\delta))+\sin(0.5(3x-2\delta))+3\delta \cos(0.5x)+6\sin(0.5(x-\delta))\big) \Big).
\end{split}
\end{equation}
Three types of distributions are considered for $\xi_{(i)}$, $i=1,\cdots,5$: the uniform distribution $\xi_{(i)}^{(1)}\sim \mathcal{U}[-0.1,0.1]$, the Gaussian distribution  $\xi_{(i)}^{(2)}\sim \mathcal{N}(0,0.1^2)$ and the lognormal distribution $\xi_{(i)}^{(3)}=\exp(\xi_{(i)}^{(2)})$. 

Numerical results are provided in Figure \ref{fig:test11,ex2it4}. With fixed Smolyak sparse grid level $\eta=5$, 
we demonstrate the spatial convergence of numerical solution for grid sizes $h=\{1/10, 1/20, 1/40, 1/80, 1/160, 1/320\}$ in Figure \ref{fig:test11,ex2it4}(a). First-order convergence $O(h)$ is observed. In Figures \ref{fig:test11,ex2it4}(b) and \ref{fig:test11,ex2it4}(c), we employ a fixed grid size $h=1/1280$ and demonstrate the convergence of solution error with increasing sparse grid level $\eta=1,\cdots,5$ in the parametric space. When $\eta\geq 3$, the numerical error of mean for uniform distribution reaches a convergence plateau because of the spatial discretization error. Noting that we have $\eta< N/\log(2)\approx 7.5$ in this case, algebraic convergence of the sparse grid PCM is verified for all three types of distributions. Therefore, the $O(h)$ spatial convergence and the algebraic convergence of the sparse grid PCM together verify the estimates in Theorem \ref{thm:nonlocalconv}.

\subsubsection*{Test 2: consistency study on a problem with 2D physical domain and 1D parametric space}

\begin{figure}[h!]
\centering
\subfigure[Convergence with refinement in the physical space.]{\includegraphics[width=.45\columnwidth]{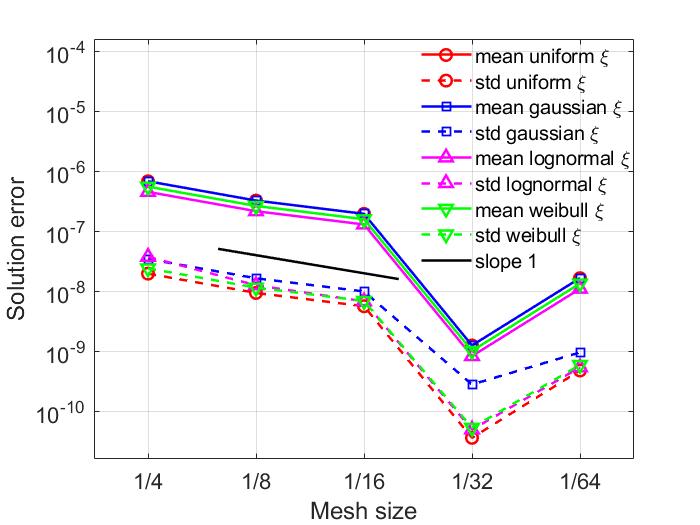}}\\
\subfigure[Convergence with sample numbers in the log scale.]{\includegraphics[width=.45\columnwidth]{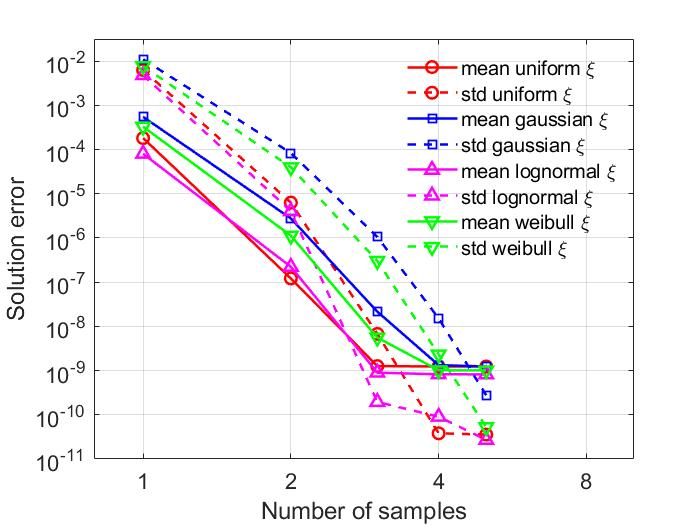}}
\subfigure[Convergence with sample numbers in the linear scale.]{\includegraphics[width=.45\columnwidth]{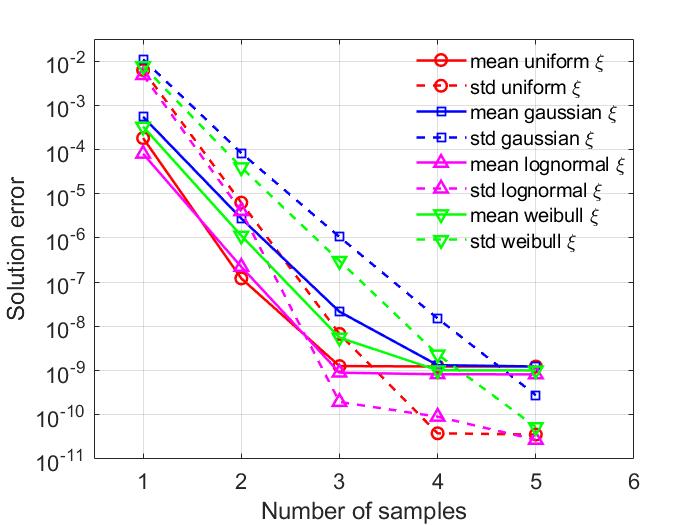}}
\caption{Test 2: consistency study of the numerical solution to the nonlocal analytical solution on a problem with 2D physical domain and 1D parametric space. Results in (a) are generated with $5$ samples, which corresponds to Smolyak formulation level $5$. The data points in (b) and (c) correspond to Smolyak formulation levels $\eta=1,\cdots,5$.}
\label{fig:test12,ex4it3}
\end{figure}


We now consider a case on a square physical domain  $\Omega=[0,1]\times[0,1]$ depending on a random variable $\xi$. With fixed $\delta=0.525$ and nonlocal diffusion coefficient $A(\xb,\yb,\xi)=(2+\xi)[2+\cos(A(x_{(1)}+y_{(1)}))\cos(B(x_{(2)}+y_{(2)}))]/\delta^4$, we consider the manufactured nonlocal analytical solution 
$$u^\delta(\xb,\xi)=u(x_{(1)},x_{(2)},\xi)=\cos(Ax_{(1)})\sin(Bx_{(2)})/(2+\xi).$$
Here we take $A=B=0.3$ in this example. Four types of popular distributions are studied: the uniform distribution $\xi^{(1)}\sim \mathcal{U}[-0.1,0.1]$, the Gaussian distribution  $\xi^{(2)}\sim \mathcal{N}(0,0.1^2)$,  the lognormal distribution {$\xi^{(3)}=\exp(\xi^{(2)})$}, and the rescaled Weibull distribution $\xi^{(4)}=0.5\hat{\xi}$. Here $\hat{\xi}$ is the Weibull random variable with the shape parameter $k=5.0$ and the scale parameter $\lambda=1$. 

Numerical convergence to the nonlocal analytical solution is demonstrated in Figure \ref{fig:test12,ex4it3}. With $5$ collocation points in PCM, we investigate the spatial convergence of numerical solution for grid sizes $h=\{1/4, 1/8, 1/16, 1/32, 1/64\}$ in Figure \ref{fig:test12,ex4it3}(a). An $O(h)$ convergence rate is observed. In Figures \ref{fig:test12,ex4it3}(b) and \ref{fig:test12,ex4it3}(c) we consider a fixed grid size $1/32$ and study the convergence of solution error with increasing sample numbers in PCM. In particular, we take sparse grid levels $\eta\in\{1,\cdots,5\}$. When $\eta\geq 4$, the numerical errors of mean for all distribution types reach a convergence plateau because of the spatial discretization error. When $\eta>N/\log(2)\approx 1.5$, sub-exponential convergence is observed for all four types of considered distributions before reaching this plateau, which again verifies Theorem \ref{thm:nonlocalconv}.

\subsection{Asymptotic Compatibility (AC) to the Stochastic Local Limit}\label{sec:localtest}

In this section we investigate the asymptotic compatibility (AC) of the proposed approach by studying the convergence of its numerical solution to the corresponding local limit when $\delta,h\rightarrow 0$. In particular, we focus on the $\delta$-convergence limit and fix the ratio between $\delta$ and $h$. Following the conventions in \cite{yu2021asymptotically,silling2021propagation,you2021data}, in all tests we consider the nonlocal diffusion coefficient as the harmonic mean of the local diffusion coefficient: $A(\xb,\yb,\xib)=2/(a^{-1}(\xb,\xib)+a^{-1}(\yb,\xib))$.

\subsubsection*{Test 1: AC study on a problem with 1D physical domain and 5D parametric space}

\begin{figure}[h!]
\centering
\subfigure[Convergence with $\delta,h\rightarrow 0$ in the physical space.]{\includegraphics[width=.45\columnwidth]{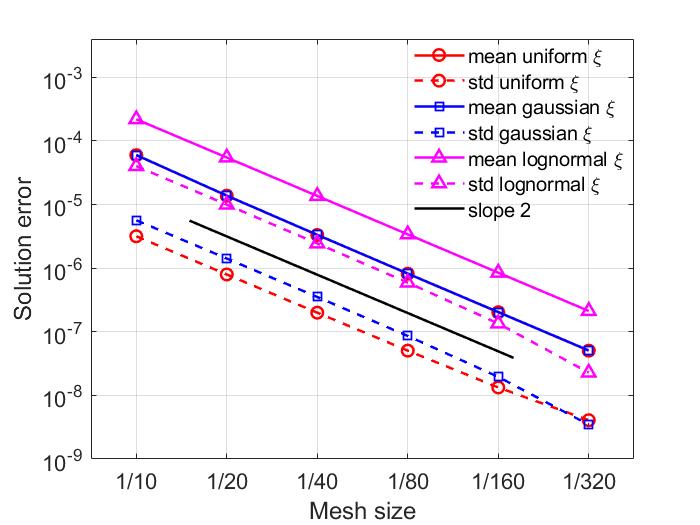}}\\
\subfigure[Convergence with sample numbers in the log scale.]{\includegraphics[width=.45\columnwidth]{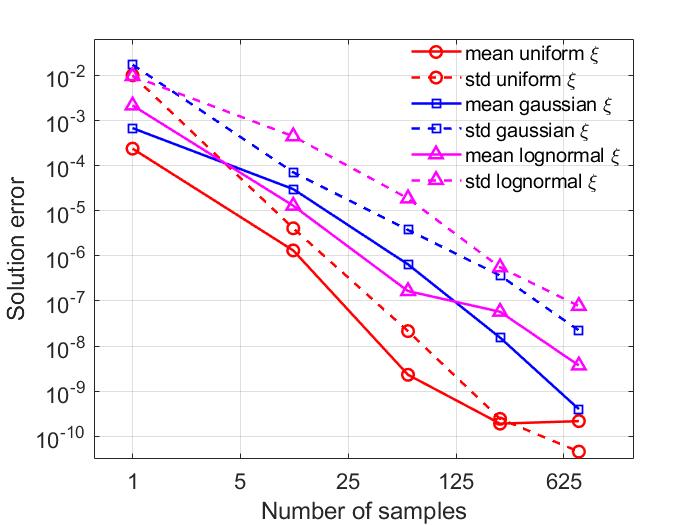}}
\subfigure[Convergence with sample numbers in the linear scale.]{\includegraphics[width=.45\columnwidth]{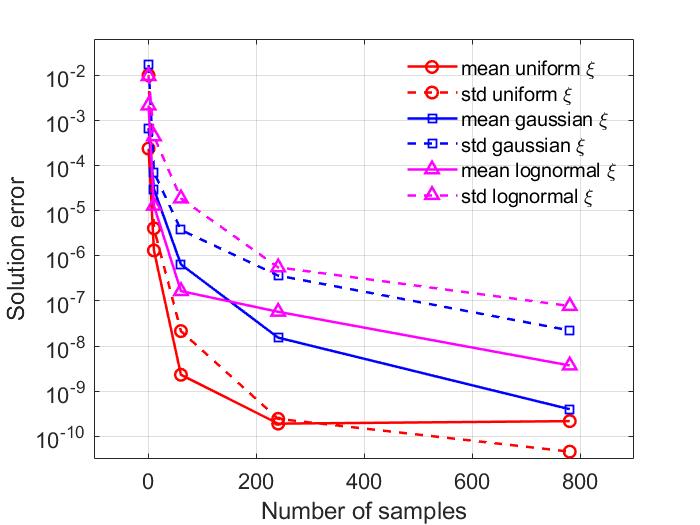}}
\caption{Test 1: asymptotic compatibility study of the numerical solution to the analytical local limit on a problem with 1D physical domain and 5D parametric space. Results in (a) are generated with $781$ samples, which corresponds to Smolyak formulation level $5$. The data points in (b) and (c) corresponds to Smolyak formulation levels $\eta=1,\cdots,5$.}
\label{fig:test21,ex3it33}
\end{figure}

We consider a case with 1D physical domain $\Omega=[-1,1]$ and 5D parametric space $\xib=(\xi_{(1)},\cdots,\xi_{(5)})$, where $\xi_{(i)}$ are i.i.d. random variables. The analytical local solution is given by
$$u^0(x,\xib)=\dfrac{\log(12+(1+\exp(\sin(\xi_{(1)}))+\cos(\xi_{(2)})+\exp(\sin(\xi_{(3)}))+\cos(\xi_{(4)})+\exp(\sin(2\xi_{(5)}))) \sin(x))}{1+\exp(\sin(\xi_{(1)}))+\cos(\xi_{(2)})+\exp(\sin(\xi_{(3)}))+\cos(\xi_{(4)})+\exp(\sin(2\xi_{(5)})))},$$
with fixed loading $f(x)=\sin(x)$ and local diffusion coefficient
$$a(x,\xib)=12+(1+\exp(\sin(\xi_{(1)}))+\cos(\xi_{(2)})+\exp(\sin(\xi_{(3)}))+\cos(\xi_{(4)})+\exp(\sin(2\xi_{(5)})))\sin(x).$$
Three types of distributions are considered for $\xi_{(i)}$, $i=1,\cdots,5$: the uniform distribution $\xi_{(i)}^{(1)}\sim \mathcal{U}[-0.1,0.1]$, the Gaussian distribution $\xi_{(i)}^{(2)}\sim \mathcal{N}(0,0.1^2)$ and the lognormal distribution $\xi_{(i)}^{(3)}=\exp(\xi_{(i)}^{(2)})$.

Numerical results are provided in Figure \ref{fig:test21,ex3it33}. With fixed ratio $\delta/h=3.8$ and Smolyak sparse grid level $\eta=5$, in Figure \ref{fig:test21,ex3it33}(a) we show the error of numerical solution with respect to the analytical local limit for grid sizes $h=\{1/10, 1/20, 1/40, 1/80, 1/160, 1/320\}$. Second-order convergence $O(\delta^2)$ is observed. The proposed approach is therefore AC and the observed convergence rate is consistent with our AC analysis in Theorem \ref{thm:localconv}. In Figures \ref{fig:test21,ex3it33}(b) and \ref{fig:test21,ex3it33}(c) we fix $h=1/5000$ and $\delta=3.8h$, and show the convergence of solution error with increasing sparse grid level $\eta=1,\cdots,5$ in the parametric space. Noting that we have $\eta< N/\log(2)\approx 7.5$ in this case, algebraic convergence of the sparse grid PCM is verified for all types of distributions, which is also consistent with the analysis in Theorem \ref{thm:localconv}.


\subsubsection*{Test 2: AC study on a problem with 2D circular domain and 1D parametric space}

\begin{figure}[h!]
\centering
\subfigure[Convergence with $\delta,h\rightarrow 0$ in the physical space.]{\includegraphics[width=.45\columnwidth]{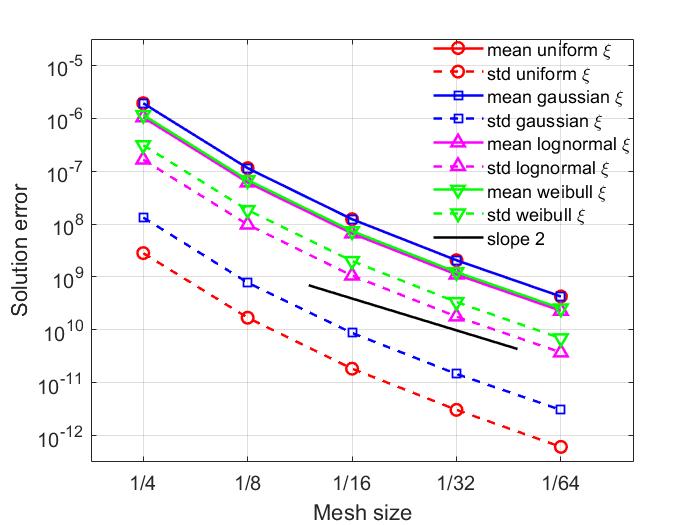}}\\
\subfigure[Convergence with sample numbers in the log scale.]{\includegraphics[width=.45\columnwidth]{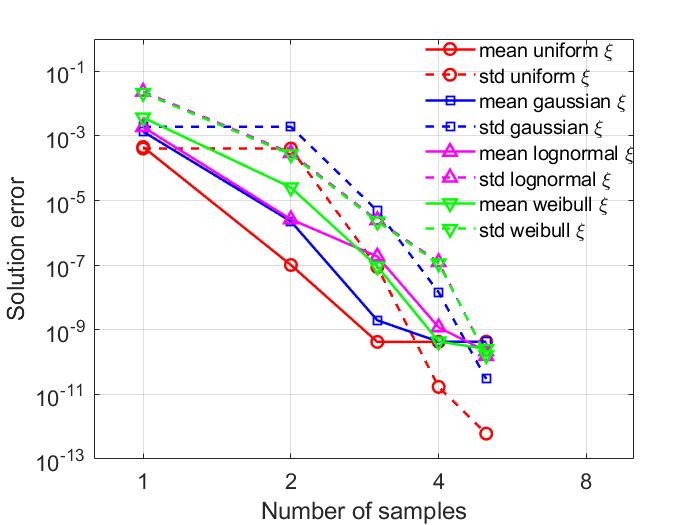}}
\subfigure[Convergence with sample numbers in the linear scale.]{\includegraphics[width=.45\columnwidth]{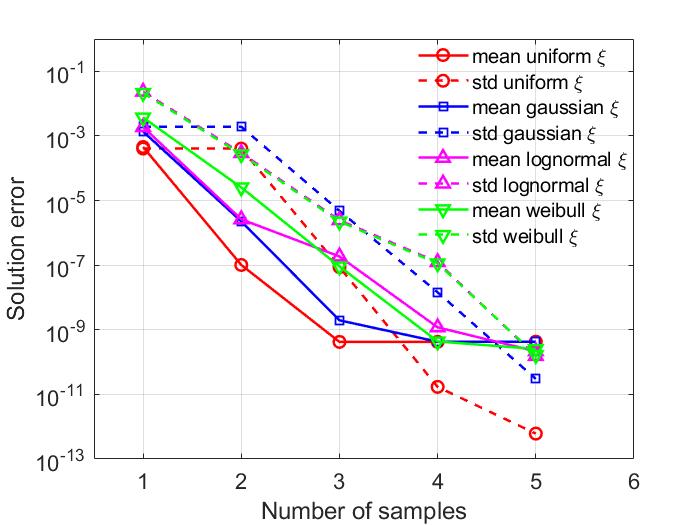}}
\caption{Test 2: asymptotic compatibility study of the numerical solution to the analytical local limit on a problem with 2D circular domain and 1D parametric space. Results in (a) are generated with $7$ samples, which corresponds to Smolyak formulation level $5$. The data points in (b) and (c) correspond to Smolyak formulation levels $\eta=1,\cdots,5$.}
\label{fig:test22,ex5it1}
\end{figure}

We next consider a more general physical domain with a curvilinear boundary. In particular, we employ a circular physical domain  $\Omega=B_1(0)$ and consider 1D parametric space. Under radical coordinate and with local diffusion coefficient $a(r,\theta,\xi)=1/(2+\cos(\xi)\sin(r^2))$, loading $f(r,\theta)=1$, the classical diffusion problem yields the analytical local solution 
$$u^0(r,\theta,\xi)=\frac{1}{4}(2r^2-\cos(\xi)\cos(r^2)).$$
Four types of distributions are considered in this case: the uniform distribution $\xi^{(1)}\sim \mathcal{U}[-0.1,0.1]$, the Gaussian distribution  $\xi^{(2)}\sim \mathcal{N}(0,0.1^2)$, the lognormal distribution $\xi^{(3)}=\exp(\xi^{(2)})$, and the Weibull distribution $\xi^{(4)}$ with the shape parameter $k=5.0$ and the scale parameter $\lambda=1$. 

Numerical results are shown in Figure \ref{fig:test22,ex5it1}. In Figure \ref{fig:test22,ex5it1}(a) we consider fixed $\delta/h = 3.8$ and $7$ samples with sparse grid PCM, then study the convergence of numerical solution to the analytical local limit with decreasing grid size from $h=1/4$ to $1/64$. An $O(\delta^2)$ is obtained. In Figures \ref{fig:test22,ex5it1}(b) and \ref{fig:test22,ex5it1}(c) we employ a fixed grid size {$h=1/64$ and $\delta=3.8h$}, then demonstrate the convergence rate of PCM with increasing number of samples corresponding to $\eta=1,\cdots,5$. Before the numerical error reaches a convergence plateau due to the spatial discretization error, an sub-exponential convergence is obtained with increasing sample numbers.

\subsubsection*{Test 3: AC study on a problem with 2D square domain and 2D parametric space}

\begin{figure}[h!]
\centering
\subfigure[Convergence with $\delta,h\rightarrow 0$ in the physical space.]{\includegraphics[width=.48\columnwidth]{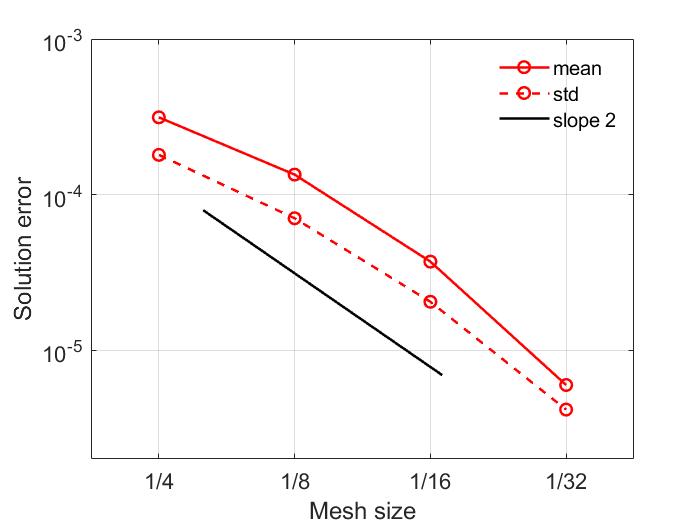}}
\subfigure[Convergence with increasing sample numbers.]{\includegraphics[width=.48\columnwidth]{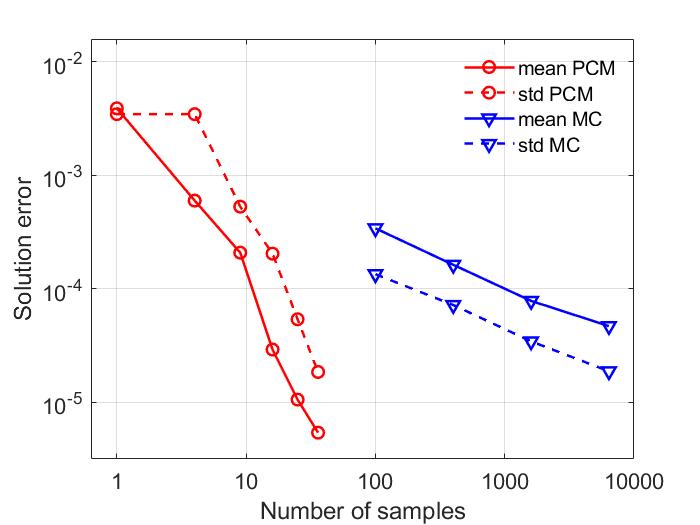}}
\caption{Test 3: asymptotic compatibility study of the numerical solution to the analytical local limit on a problem with 2D square domain and 2D parametric space. Results in (a) are generated with $100$ samples. The data points in (b) correspond to Smolyak formulation levels $\eta=1,\cdots,6$.}
\label{fig:ex5it2}
\end{figure}

We now consider a 2D square domain $\Omega=[-1, 1]\times [-1, 1]$ and 2D parametric space $\xib=(\xi_{(1)},\xi_{(2)})$, $\xi_{(i)}\sim \mathcal{U}[-0.1,0.1]$, with loading $f(\xb)=1$ and local diffusion coefficient 
$$a(\xb,\xib)=3+\sum_{k=1}^2\frac{\cos(30 \xi_{(k)})-1}{k^2}\cos(2kx_{(1)})\sin(2kx_{(2)}).$$
Homogeneous Dirichlet-type boundary condition is applied on boundary $\partial \Omega$. Since there is no analytical expression for the local limit, $u^0(\xb,\xib)$ is generated numerically based on a spectral method solver. The mean and standard deviation of the analytical local solution is calculated with the Monte Carlo (MC) method with $100,000$ samples.  

Numerical results are provided in Figure \ref{fig:ex5it2}. With $100$ collocation points in PCM and fixed $\delta/h = 2.8$, we investigate the spatial convergence of numerical solution with decreasing grid size from $h=1/4$ to $1/32$ in Figure \ref{fig:ex5it2}(a). Second order convergence is observed,  which is consistent with our analysis. In Figure \ref{fig:ex5it2}(b) we employ fixed grid size $h=1/32$, $\delta=2.8h$ and show the convergence of numerical solutions obtained from MC simulations and the PCM results with respect to the increase of sample points. Here the sparse grid levels are taken as $\eta=1,\cdots,6$. When $\eta>N/\log(2)\approx3$, sub-exponential convergence is observed for PCM, and the results also indicate that PCM can achieve a similar accuracy with much smaller number of sample points than MC.

\section{Stochastic Nonlocal Diffusion Problem in Randomly Heterogeneous Domain}\label{sec:exp}

\subsection{Stochastic Representation and Karhunen-Lo\`eve Expansion}\label{sec:KL}

We use the Karhunen-Lo\`eve (K-L) expansion to represent the random field $a(\bm x, \omega)$. In general, consider a square-integrable stochastic random field $F(\bm x, \omega)$ defined on $\Omega\times \Omega_p$, where $\Omega$ is a subset of $\mathbb{R}^d$ ($d$ is the dimension) and $\Omega_p$ is the sample space of a probability space $(\Omega_p,\mathcal{F},\mathcal{P})$. If $F(\bm x, \omega)$ has a constant mean and a continuous covariance function (also called kernel function) $\Xi(\bm x, \bm y)$, then $F(\bm x, \omega)$ can be represented by the following Karhunen-Lo\`eve (KL) expansion:
\[F(\bm x, \omega)=F_0+\sum_{n=1}^{\infty}\sqrt{\lambda_{i}}\phi_{i}(\bm x)\xi_{(i)}(\omega),\]
where $F_0$ is the constant mean, $(\lambda_{i},\phi_{i})$ are eigenpairs (i.e., eigenvalue and corresponding eigenfunction) of the kernel function $\Xi$, and $\xi_{(i)}$ are independent random variables with zero mean and unit variance. In practice, the summation is truncated up to $N$ terms for computational purpose, where $N$ is taken such that
\begin{equation}\label{eqn:trunc}
\sum_{i=1}^N\lambda_i \geq 0.9 \sum_{i=1}^\infty\lambda_i.
\end{equation}
Notice that instead of $0.9$, other numbers like $0.85$ and $0.95$ are also widely used. 
The notation of $F(\bm x,\omega)$ is then replaced with $F(\bm x, \bm\xi)$, where $\bm\xi=(\xi_{(1)}, \xi_{(2)},\dotsc,\xi_{(N)})$.
Here, we use the KL expansion to represent the random field $a(\bm x, \bm\xi)$, and we consider the nonlocal diffusion coefficient as the harmonic mean of the local diffusion coefficient: $A(\xb,\yb,\xib)=2/(a^{-1}(\xb,\xib)+a^{-1}(\yb,\xib))$ in our nonlocal model.

\subsection{Numerical Simulations}

\begin{figure}[h!]
\centering
\subfigure[Convergence with $\delta,h\rightarrow 0$ in the physical space.]{\includegraphics[width=.48\columnwidth]{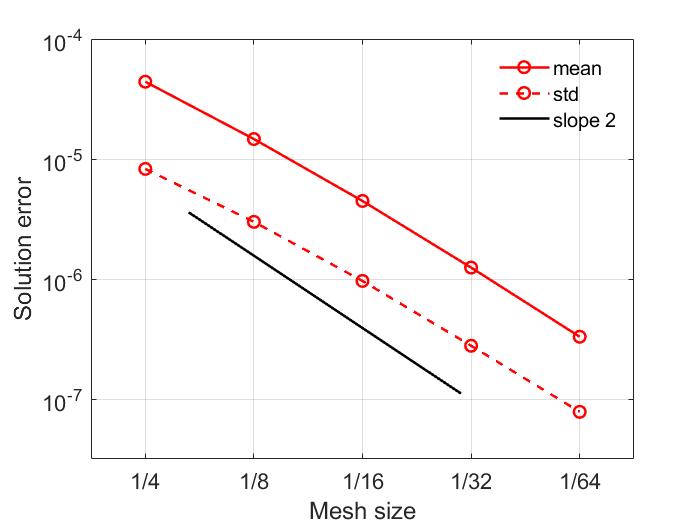}}
\subfigure[Convergence with increasing sample numbers.]{\includegraphics[width=.48\columnwidth]{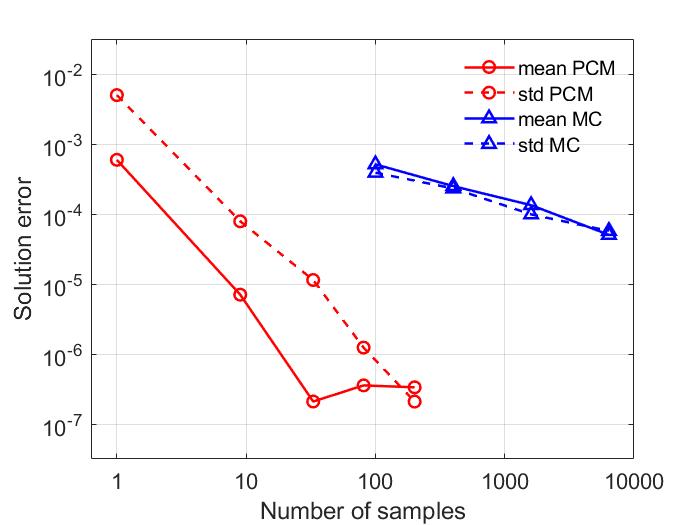}}
\caption{Asymptotic compatibility study of the numerical solution to the analytical local limit on randomly heterogeneous nonlocal problem with a given  spatial correlation structure. Results in (a) are generated with $441$ samples, which corresponds to Smolyak formulation level $5$. The data points in (b) correspond to Smolyak formulation levels $\eta=1,\cdots,5$.}
\label{fig:ex6it1}
\end{figure}

In this section, we consider the proposed approach on a randomly heterogeneous nonlocal problem with a given  spatial correlation structure. We consider a 2D domain $\Omega=[-1,\, 1]\times [-1, \, 1]$ with $u^0=0$ on $\partial \Omega$. The local diffusivity coefficient $a(\xb,\omega)$ is modeled as a random field with constant mean $a_0=4$ and a continuous covariance function:
$$\Xi(\xb, \yb)(=\text{Cov}(a(\xb,\xib),a(\yb,\xib)))=\sigma^2\exp\left(-\frac{|x_{(1)}-y_{(1)}|^2}{\eta_1}-\frac{|x_{(2)}-y_{(2)}|^2}{\eta_2}\right),$$
with $\sigma=1,\eta_1=\eta_2=1$. Of note, $\sigma^2$ is the variance and $\eta_i, i=1,2$ are correlation lengths. We note that the above covariance kernel is separable, and therefore the eigenvalues and corresponding eigenfunctions in $\Omega$ can be derived by the multiplicaitons of the eigenvalues and eigenfunctions in the one-dimensional case. In particular, we first considerthe eigenvalues and corresponding eigenfunctions in K–L expansion for the 1D covariance function $\Xi^{(k)}(x, y)=\sigma\exp\left(-\frac{|x-y|^2}{\eta_k}\right)$, $k=1,2$ and obtain the 1D eigenpairs $(\lambda^{(k)}_{i},\phi^{(k)}_{i})$. The random local diffusivity coefficent is then expressed as the following expansion:
$$a(\xb,\bm\xi)=4+\sum_{i=1}^{N^{(1)}}\sum_{j=1}^{N^{(2)}} \sqrt{\lambda_{i}^{(1)}}\sqrt{\lambda_{j}^{(2)}}\phi^{(1)}_{i}(x_{(1)})\phi_{j}^{(2)}(x_{(2)})\xi_{(i,j)}.$$
In this computational example, $N^{(1)}=N^{(2)}=2$ is required to achieve the truncation error criteria \eqref{eqn:trunc}. We assume $\xi_{(i,j)}$, $i,j=1,2$, to be Gaussian random variables: $\xi_{(i,j)} \sim N(0,1)$, then investigate the convergence of $u^\delta_{h,K}$ to the local limit. Note that in this problem there exists no analytical expression for the local solution, we therefore generate $u^0(\xb,\xib)$ numerically based on a spectral method solver. Since the MC method suffers from a very slow convergence and requires more than $10^7$ samples to achieve a desired precision, the mean and standard deviation of the local limit is calculated with  the sparse grid PCM method with $2881$ samples (sparse grid level $\eta=10$).

Numerical results are provided in Figure \ref{fig:ex6it1}. In Figure \ref{fig:ex6it1}(a) we demonstrate the convergence in the physical space with uniform grids $h\in\{1/4, 1/8, 1/16, 1/32,1/64\}$ and fixed $\delta/h = 2.8$, using $441$ sample points (sparse grid level $\eta=5$). Second order convergence is observed, which verifies the analysis in Theorem \ref{thm:localconv}. To demonstrate the convergence in the parametric space with increasing sample numbers, in Figure \ref{fig:ex6it1}(b) we employ a fixed grid size $h=1/64$, $\delta=2.8h$, and show the numerical error to the local limit with sparse grid levels $\eta=1,\cdots,5$ in PCM. The results again indicate that comparing with MC, the sparse grid PCM achieves a better accuracy with far smaller number of samples. An algebraic convergence is observed which is consistent with Theorem \ref{thm:localconv}.

\section{Summary and Discussion}\label{sec:conclusion}


Due to the limitation of computational resources and experimental resolutions, in many applications the physical parameters describing continuous properties of heterogeneous materials cannot be accurately characterized in all details. This problem becomes more acute in the nonlocal setting, due to its relatively lack of sparsity and correspondingly larger computational cost. 

In this work we aim to consider the spatial variability of material properties in nonlocal models by solving randomly heterogeneous nonlocal diffusion problems. In particular, we have proposed an asymptotically compatible stochastic numerical method for randomly heterogeneous nonlocal diffusion problems, and provided rigorous mathematical analysis and error estimates. For spatial discretization, a meshfree discretization method with optimization-based quadrature rule is employed, which presents an up to $O(h)$ consistency error to the nonlocal solution and an $O(\delta^2)$ convergence to the local limit. On the random {parametric} space, a probabilistic collocation method (PCM) with sparse grids is employed to sample the stochastic process. Since the fast convergence of the sparse grid PCM approach relies on smoothness of the solution in the random {parametric} space, we have proved that the nonlocal solution is analytic in the input random variable and therefore guarantees an at least algebraic convergence with increasing sample numbers. This work has for the first time provided a rigorous and comprehensive mathematical framework to add uncertainty quantification analysis onto pre-existing deterministic codes for nonlocal problems with guaranteed convergence in both physical and parametric spaces.

Future directions of this research will include the development and analysis of a generalization of this approach to the nonlocal mechanics problems, such as the peridynamics. A further extension of the meshfree discretization approach will also be pursued, to achieve a higher order convergence rate.

\section*{Acknowledgements}

 Y. Fan and Y. Yu would like to acknowledge support by the National Science Foundation under award DMS 1753031 and Lehigh's High Performance Computing systems for providing computational resources at Sol. Portions of this research were conducted on Lehigh University's Research Computing infrastructure partially supported by NSF Award 2019035. X. Tian's research is support in part by the National Science Foundation grant DMS-2111608. X. Li's research is support in part by NSF DMS-1720245. X. Yang's research is support in part by the Energy Storage Materials Initiative, which is a Laboratory Directed Research and Development Project at Pacific Northwest National Laboratory.

\bibliographystyle{elsarticle-num}
\bibliography{yyu}

\end{document}